\documentclass[12pt]{article}

\usepackage{a4wide,epsfig}
\usepackage{amsfonts,amsmath}
\usepackage{caption}
\usepackage[labelformat=simple]{subcaption}

\usepackage{enumerate}
\usepackage{pgfplotstable}
\usepackage{booktabs}

\newtheorem{theorem}{Theorem}[section]

\newtheorem{lemma}[theorem]{Lemma}

\newtheorem{remark}{Remark}[section]
\newcommand{\proof} [1]
   { \noindent {\bf Proof.} #1 \hfill\rule{0.5em}{1.2ex} \par\medskip}

\newcommand{\R}{\mathbb{R}}

\newcommand{\spf}[2]{{\left\langle{#1},{#2}\right\rangle}}
\newcommand{\norm}[1]{{\left\lVert{#1}\right\rVert}}

\numberwithin{equation}{section} 

\begin{document}

\setcounter{page}{1}

\title{Adaptive least-squares space-time \\
  finite element methods}
\author{Christian~K\"othe, Richard~L\"oscher, Olaf~Steinbach}
\date{Institut f\"ur Angewandte Mathematik, TU Graz, \\[1mm]
  Steyrergasse 30, 8010 Graz, Austria \\[1mm]
  {\footnotesize \tt c.koethe@tugraz.at,
    loescher@math.tugraz.at, o.steinbach@tugraz.at}
}

\maketitle

\begin{abstract}
  We consider the numerical solution of an abstract operator equation $Bu=f$ by
  using a least-squares approach. We assume that $B: X \to Y^*$ is
  an isomorphism, and that $A : Y \to Y^*$ implies a norm in $Y$, where
  $X$ and $Y$ are Hilbert spaces. The minimizer of the least-squares
  functional $\frac{1}{2} \, \| Bu-f \|_{A^{-1}}^2$,
  i.e., the solution of the operator equation, is then
  characterized by the gradient equation $Su=B^* A^{-1}f$
  with an elliptic and
  self-adjoint operator $S:=B^* A^{-1} B : X \to X^*$. When introducing
  the adjoint $p = A^{-1}(f-Bu)$ we end up with a saddle point
  formulation to be solved numerically by using a mixed finite element
  method. Based on a discrete inf-sup stability condition we derive
  related a priori error estimates. While the adjoint $p$ is zero by
  construction, its approximation $p_h$ serves as a posteriori error
  indicator to drive an adaptive scheme when discretized appropriately.
  While this approach can be applied to rather general equations, here we
  consider second order linear partial differential equations,
  including the Poisson equation, the heat equation, and the wave equation,
  in order to demonstrate its potential, which
  allows to use almost arbitrary space-time finite element methods
  for the adaptive solution of time-dependent partial differential equations.
\end{abstract}


\section{Introduction}
The use of Galerkin finite element methods for the numerical solution
of partial differential equations is well established, see, e.g., the
text books \cite{BrennerScott:2008,Steinbach:2008}, and many others.
However, for time dependent problems it is a common procedure to first
discretize the spatial part using a finite element method, and then
applying either a time stepping method, e.g., an explicit or implicit Euler
scheme, see, e.g., the monograph \cite{Thomee:2006}, or
discontinuous Galerkin methods in time, see, e.g.,
\cite{Neumueller:2013}, and the references given therein. Recently,
the interest of discretizing space
and time at once has been rising, resulting in so called space-time
discretization methods, see, e.g., \cite{LangerSteinbach:2019}.
Although, the discretization in space and time leads to
larger systems of algebraic equations to be solved, these methods bring
the advantage of having full control of the discretization in space and time
simultaneously, allowing for space-time adaptivity. Moreover,
space-time methods offer more flexibility in the construction of efficient
solvers than time stepping methods, since preconditioning and parallelization 
in the space-time domain is applicable, and mandatory, see, e.g.,
\cite{GanderNeumueller:2016, LangerZank:2021}
for space-time solvers in the case of parabolic equations.
The derivation of
space-time formulations usually results in Petrov-Galerkin schemes, see, e.g.,
\cite{Andreev:2013, SchwabStevenson:2009, Steinbach:2015, UrbanPatera:2014}
in the case of the heat equation, where for the numerical treatment
it is crucial to establish a
related discrete inf-sup stability condition. This becomes even more
involved in the case of the wave equation where a CFL condition 
is required, e.g., \cite{SteinbachZank:2020}. Possible approaches to
overcome such a restriction is the use of discontinuous Galerkin methods, e.g.,
\cite{DoerflerWienersZiegler:2020,Neumueller:2013, NeumuellerSteinbach:2011},
or introducing a suitable transformation operator such as the
(modified) Hilbert transformation
\cite{LoescherSteinbachZank:2022,PerugiaSchwabZank:2023,SteinbachZank:2020}.
Another approach is to replace the direct variational formulation by a
least-squares/minimal residual equation. This has been studied in the context
of first order least-squares systems in, e.g., \cite{Fuehrer:2022,
  FuehrerKarkulik:2021, Fuehrer:2023, SchafelnerVassilevski:2021},
and in the context of minimal residual Petrov--Galerkin discretizations in 
\cite{Andreev:2013,StevensonWesterdiep:2021}, just to mention a few. For the
latter approach it is well known that the Galerkin discretization results
in a mixed system, where the second variable is the Riesz lift of the
residual of the primal variable. This is also the point of view we will take.

From the well-established theory of least-squares methods, with a
tremendous overview by Bochev and Gunzberger \cite{BochevGunzburger:2009},
see also \cite{BochevGunzburger:2006} and 
\cite{BochevDemkowiczGopalakrishnanGunzburger:2014,BochevGunzburger:2016},
least-squares formulations of second order partial differential
equations come with the advantage of offering an error estimator for free,
but also double the degrees of freedom. To apply the theory, it is of main
interest to consider so called practical methods, that allow to measure
the residual in a localizable norm, i.e., in a non-negative Sobolev norm
allowing at least $L^2$ regularity. Trimmed to this interest,
Führer and Karkulik introduced a first order system least-squares
method (FOSLS) for parabolic problems \cite{FuehrerKarkulik:2021}, showing
the applicability and power of least-squares methods in the space-time
setting. However, the reformulation as a first order system comes with the 
fact that one has to assume a higher regularity on the source term in order 
to show convergence of the method. This is unnatural to a certain extent, 
since
source terms of partial differential equations are usually
considered as functionals acting on the test space and thus belong to
Sobolev spaces of negative order.
%Though, source terms of partial differential equations are usually
%considered as functionals acting on the test space, and thus belong to
%Sobolev spaces of negative order, to be able to apply a practical
%least-squares framework one has to reformulate the partial differential
%equation as a first order system with the disadvantages of introducing
%further unknowns and the necessity of assuming higher regularity on the
%source term in order to show convergence of the method.
At this point we
want to mention that recently in \cite{Fuehrer:2022,Fuehrer:2023},
FOSLS with a source term in a negative order Sobolev space
were analyzed by replacing the load by its finite element approximation.  

In this paper, we will formulate and analyze a least-squares method for
the solution of abstract operator equations, that overcomes this problem,
i.e., we will be able to phrase the method, even when the source is of
minimal regularity. For $X$ and $Y$ being Sobolev spaces of non-negative
order, we will consider the problem to find $u\in X$ such that $Bu=f \in Y^*$,
where $B : X \to Y^*$ is an isomorphism. Then the problem is equivalent to
minimize the residual $\frac{1}{2} \| Bv-f \|_{Y^*}^2$ over all $v\in X$.
As $Y^*$ is the dual space of $Y$, the residual is measured in a Sobolev
norm with negative index. The main idea is now based on lifting the
negative order Sobolev norm using a bijective operator $A:Y\to Y^*$
in order to keep the method practical. Then the problem is equivalent
to minimize the functional $\frac{1}{2} \| Bv-f \|_{A^{-1}}^2$, for which
the minimizer $u\in X$ is characterized as the unique solution of the
gradient equation $B^* A^{-1}(Bu-f) = 0$. Introducing the auxiliary
variable $p=A^{-1}(f-Bu) \in Y$ we end up with solving a saddle point
formulation. This saddle point point problem is of similar
shape as those obtained from optimal control with energy regularization, see 
\cite{LSTY:2021,2022_Langer_Steinbach_robust,LoescherSteinbach:2022}.
The discretization follows by standard means, using conforming
finite element spaces of lowest order. To ensure uniqueness, we will
need a discrete inf-sup stability condition. Though, $p\equiv 0$ on the
continuous level, for the discrete Lagrange multiplier in general it holds
that $p_h\neq 0$. Thus, it can be used to define an a posteriori error
estimator which is shown to be efficient and reliable using a so called
saturation assumption 
\cite{2016_Carstensen_JustificationOfTheSaturationAssumption,
  2002_Doerfler_SmallDataOscillation}.
The idea of lifting the norm was already considerd in
\cite{BrambleLazarovPasciak:1997,BrambleLazarovPasciak:1998} for elliptic
partial differential equations, or more recently in
\cite{MonsuurStevensonStorn:2023} for a larger class of problems,
and in \cite{DahmenMonsuurStevenson:2022} adding a stabilization term
when considering ill-posed problems. In \cite{Andreev:2013}, this concept
was already applied in the context of space-time methods to parabolic
equations, with a discretization on a tensor product mesh.        
The novelty of our method is the application to space-time problems on
completely unstructured meshes, and the practicability even when using
negative Sobolev norms as well as the proof of an efficient and reliable
error estimator in this setting. Moreover, the use of completely
unstructured space-time meshes allows for full space-time adaptivity,
and later on, for a parallel solution of the resulting linear systems
of algebraic equations.

This paper is organized as follows. In Section 2 we consider the
least-squares approach for the solution of an operator equation $Bu=f$
in an abstract sense. We derive the saddle point formulation and show
ellipticity of the related operator $S=B^*A^{-1}B$. Based on a discrete
inf-sup  condition we show unique solvability of the discrete system
and derive related a priori error estimates. Using a so called saturation
assumption we are able to prove that the use of the dual $p_h$ as
error estimator is efficient and reliable.
In Section 3 we apply the approach to the Poisson equation,
showing that it is also related to the well known
$h-\frac{h}{2}$ estimator \cite{Praetorius:2014,
  2010_FerrazOrtnerPraetorius_ConvergenceofSimpleAdaptiveGalerkin}.
Section~4 deals with the application of the approach to the heat equation,
while the wave equation is considered in Section~5. In all cases
we provide numerical examples to illustrate our theoretical
findings. Finally, in Section 6 we give a conclusion, where
we mention possible extensions and further work which needs to be done.
     
\section{Abstract setting}\label{sec:abstract setting}
Let $X \subset H \subset X^*$ and $Y \subset H \subset Y^*$ be Gelfand
triples of Hilbert spaces, where $X^*$, $Y^*$ are the duals of $X$, $Y$
with respect to $H$ and with the duality pairing $\langle f , q \rangle_H$
for $f \in Y^*$ and $q \in Y$. Let $A : Y \to Y^*$ be a bounded linear
operator which
is assumed to be self-adjoint and elliptic in $Y$. Therefore, $A$ implies
a norm in $Y$, i.e.,
$\| q \|_Y := \sqrt{\langle A q , q \rangle_H}$ for $q \in Y$.
For $ f \in Y^*$, the norm is given by duality,
\begin{equation}\label{Abstract dual norm}
  \| f \|_{Y^*} := \sup\limits_{0 \neq q \in Y}
  \frac{\langle f , q \rangle_H }{\| q \|_Y} \, .
\end{equation}

\begin{lemma}
The dual norm \eqref{Abstract dual norm} allows the representations
\begin{equation}\label{Abstract representation dual norm}
  \| f \|_{Y^*}^2 = \| p_f \|_Y^2 =
  \langle A p_f , p_f \rangle_H = \langle f , p_f \rangle_H =
  \langle A^{-1} f , f \rangle_H ,
\end{equation}
where $ p_f \in Y$ is the unique solution of the variational formulation
\begin{equation}\label{Abstract definition pf}
  \langle A p_f , q \rangle_H = \langle f , q \rangle_H \quad
  \mbox{for all} \; q \in Y .
\end{equation}
\end{lemma}
\proof{Since $A : Y \to Y^*$ is bounded and elliptic,
  $p_f = A^{-1} f$ is well defined as unique solution of the
  variational formulation \eqref{Abstract definition pf}.
  From the definition \eqref{Abstract dual norm}, we first have,
  note that $\| f \|_{Y^*} > 0$ implies $\| p_f \|_Y >0$,
  \[
    \| f \|_{Y^*} = \sup\limits_{0 \neq q \in Y}
    \frac{\langle f , q \rangle_H }{\| q \|_Y} \geq
    \frac{\langle f , p_f \rangle_H}{\| p_f \|_Y}, \quad
  \mbox{i.e.}, \quad
    \langle f , p_f \rangle_H \leq \| f \|_{Y^*} \| p_f \|_Y \, .
  \]
  Hence we obtain
  \[
    \| p_f \|_Y^2 = \langle A p_f , p_f \rangle_H =
    \langle f , p_f \rangle_H \leq \| f \|_{Y^*} \| p_f \|_Y, \quad
    \mbox{i.e.}, \quad
    \| p_f \|_Y \leq \| f \|_{Y^*} \, .
  \]
  On the other hand, \eqref{Abstract dual norm} implies, when
  using \eqref{Abstract definition pf},
    \[
    \| f \|_{Y^*} = \sup\limits_{0 \neq q \in Y}
    \frac{\langle f , q \rangle_H }{\| q \|_Y} =
    \sup\limits_{0 \neq q \in Y}
    \frac{\langle A p_f , q \rangle_H}{\| q \|_Y} \leq \| p_f \|_Y \, ,
  \]
  and hence, $ \| f \|_{Y^*} = \| p_f \|_Y $
  follows. With this we finally obtain
  \[
    \| f \|_{Y^*}^2 = \| p_f \|^2_Y =
    \langle A p_f , p_f \rangle_H = \langle f , p_f \rangle_H =
    \langle f , A^{-1} f \rangle_H .
  \]
}

\noindent
Let $B:X \to Y^*$ be a bounded
linear operator which satisfies an inf-sup condition, i.e.,
there exist positive constants $c_1^B$ and $c_2^B$ such that
\begin{equation}\label{Assumptions B}
  \| B v \|_{Y^*} \leq c_2^B \, \| v \|_X, \quad
  \sup\limits_{0 \neq q \in Y} \frac{\langle Bv,q \rangle_H}{\| q \|_Y}
  \geq c_1^B \, \| v \|_X \quad \mbox{for all} \; v \in X .
\end{equation}
In addition, we assume that $B$ is surjective. Then, 
$B:X \to Y^*$ is an isomorphism.
Due to the assumptions made we conclude
unique solvability of the operator equation to find $u \in X$ such
that $Bu=f$ in $Y^*$ is satisfied. The solution of the operator
equation $Bu=f$ in $Y^*$ is equivalent to the minimization of a quadratic
functional for $v \in X$,
\begin{eqnarray*}
  {\mathcal{J}}(v)
  & = & \frac{1}{2} \, \| B v - f \|_{Y^*}^2 \, = \,
        \frac{1}{2} \, \langle A^{-1} (Bv-f) , Bv-f \rangle_H \\
  & = & \frac{1}{2} \, \langle B^* A^{-1} B v,v \rangle_H
        - \langle B^* A^{-1} f , v \rangle_H +
        \frac{1}{2} \, \langle A^{-1} f,f\rangle_H,
\end{eqnarray*}
whose minimizer $u \in X$ is given as solution of the gradient equation
\[
B^* A^{-1} (Bu-f) = 0 \, ,
\]
i.e., we have to solve the operator equation
\begin{equation}\label{Operator equation Su=g}
Su := B^* A^{-1} B u = B^* A^{-1} f \, .
\end{equation}
Note that $B^* : Y \to X^*$ is the adjoint of $B:X \to Y^*$, i.e.,
\[
  \langle B^* q , v \rangle_H :=
  \langle q , B v \rangle_H \quad \mbox{for all} \; v \in X, \; q \in Y ,
\]
satisfying
\[
  \| B^* q \|_{X^*} = \sup\limits_{0 \neq v \in X}
  \frac{\langle B^* q , v \rangle_H}{\|v\|_X}
  = \sup\limits_{0 \neq v \in X}
  \frac{\langle q , B v \rangle_H}{\|v\|_X} \leq
  \sup\limits_{0 \neq v \in X}
  \frac{\| q \|_Y \| B v \|_{Y^*}}{\|v\|_X} \leq c_2^B \, \| q \|_Y 
\]
for all $q \in Y$.

\begin{lemma}\label{Lemma S}
  The operator $S:=B^* A^{-1} B : X \to X^*$ is bounded and elliptic, i.e.,
  \[
    \| S u \|_{X^*} \leq c_2^S \, \| u \|_X , \quad
    \langle S u , u \rangle_H \geq c_1^S \, \| u \|_X^2 \quad
    \mbox{for all} \; u \in X,
  \]
  where $c_2^S = [c_2^B]^2$, $c_1^S = [c_1^B]^2$.
\end{lemma}
\proof{For $u \in X$ we first have, using
  \eqref{Abstract representation dual norm} for $f = Bu$,
  \begin{eqnarray*}
    \| S u \|_{X^*}
    & = & \sup\limits_{0 \neq v \in X} \frac{\langle S u , v \rangle_H}
          {\| v \|_X} =
          \sup\limits_{0 \neq v \in X}
          \frac{\langle A^{-1} B u , B v \rangle_H}{\| v \|_X} \\
    &  \leq & \sup\limits_{0 \neq v \in X}
              \frac{\| A^{-1} B u \|_Y \| B v \|_{Y^*}}{\| v \|_X} \, \leq \,
              \sup\limits_{0 \neq v \in X}
              \frac{\| B u \|_{Y^*} \| B v \|_{Y^*}}{\| v \|_X} \, \leq \,
              [c_2^B]^2 \, \| u \|_X \, .
  \end{eqnarray*}
  In addition, we define $p_u = A^{-1} B u \in Y$ to obtain
  \[
    \| p_u \|_Y^2 = \langle A p_u , p_u \rangle_H =
    \langle B^* A^{-1} B u , u \rangle_H = \langle S u , u \rangle_H .
  \]
  From the inf-sup stability condition in \eqref{Assumptions B}
  we then conclude
  \[
    c_1^B \, \| u \|_X \leq \sup\limits_{0 \neq v \in Y}
    \frac{\langle B u , v \rangle_H}{\| v \|_Y} =
    \sup\limits_{0 \neq v \in Y} \frac{\langle A p_u , v  \rangle_H}
    {\| v \|_Y} \leq \| p_u \|_Y ,
  \]
  i.e.,
  \[
    [c_1^B]^2 \, \| u \|_X^2  \leq \| p_u \|_Y^2 \, = \,
    \langle S u ,u \rangle_H .
  \]
}

\noindent
In fact, $S := B^* A^{-1} B : X \to X^*$ defines an equivalent norm in $X$,
\[
  \| v \|_S := \sqrt{\langle S v , v \rangle_H} =
  \sqrt{\langle A^{-1} B v , B v \rangle_H} = \| B v \|_{Y^*} \, ,
\]
satisfying
\begin{equation}\label{norm equivalence inequalities}
  c_1^B \, \| v \|_X \leq \| v \|_S \leq c_2^B \, \| v \|_X \quad
  \mbox{for all} \; v \in X .
\end{equation}
The variational formulation of the operator equation
\eqref{Operator equation Su=g} is to find $u \in X$ such that
\begin{equation}\label{VF Operator equation Su=g}
  \langle S u , v \rangle_H = \langle B^* A^{-1} f , v \rangle_H
  \quad \mbox{for all} \; v \in X ,
\end{equation}
which is uniquely solvable for all $f \in Y^*$.

Let $X_H = \mbox{span} \{ \varphi_k \}_{k=1}^{M_X} \subset X$ be some finite
dimensional ansatz space which is defined with respect to some
admissible decomposition of the computational domain
into shape-regular simplicial finite elements of mesh size $H$,
see, e.g., \cite{BrennerScott:2008, Steinbach:2008}.
Then the Galerkin variational formulation of
\eqref{VF Operator equation Su=g} is to find $u_H \in X_H$ such that
\begin{equation}\label{Galerkin VF Operator equation Su=g}
  \langle S u_H , v_H \rangle_H = \langle B^* A^{-1} f , v_H \rangle_H
  \quad \mbox{for all} \; v_H \in X_H.
\end{equation}
Using standard arguments, we conclude Cea's lemma,
\begin{equation}\label{Cea Operator equation Su=g}
  \| u - u_H \|_S \leq \inf\limits_{v_H \in X_H} \| u - v_H \|_S \, ,
\end{equation}
and convergence $u_H \to u$ in $X$ follows from an approximation
property of $X_H$.

Since the operator $S = B^* A^{-1} B$ does not allow a direct evaluation
in general, we have to define a suitable approximation.
For $u \in X$ we write $S u = B^* A^{-1} B u = B^* p_u$, where
$ p_u = A^{-1} B u$. In fact, $p_u \in Y$ is the unique solution of
the variational formulation
\[
  \langle A p_u , q \rangle_H = \langle B u , q \rangle_H \quad
  \mbox{for all} \; q \in Y .
\]
Let $Y_h = \mbox{span} \{ \psi_i \}_{i=1}^{M_Y} \subset Y$ be a finite
dimensional ansatz space which is
defined with respect to some mutually different decomposition of the
computational domain into finite elements of mesh size $h$.
Then we define $p_{uh} \in Y_h$ as unique solution of the
Galerkin variational formulation
\[
  \langle A p_{uh} , q_h \rangle_H = \langle B u , q_h \rangle_H \quad
  \mbox{for all} \; q_h \in Y_h ,
\]
and we define $\widetilde{S} u := B^* p_{uh}$ as approximation of
$Su=B^* p_u$. We immediately have the bounds
\[
  \| p_{uh} \|_Y^2 = \langle A p_{uh} , p_{uh} \rangle_H =
  \langle B u , p_{uh} \rangle_H \leq \| B u \|_{Y^*} \| p_{uh} \|_Y \leq
  c_2^B \, \| u \|_X \| p_{uh} \|_Y,
\]
i.e.,
\[
\| p_{uh} \|_Y \leq c_2^B \, \| u \|_X ,
\]
and
\begin{equation}\label{Stilde bounded}
  \| \widetilde{S} u \|_{X^*} =
  \| B^* p_{uh} \|_{X^*} \leq c_2^B \, \| p_{uh} \|_Y \leq
  [c_2^B]^2 \, \| u \|_X,
\end{equation}
as well as the error estimate
\begin{equation}\label{Error Stilde}
  \| (S-\widetilde{S})u \|_{X^*} =
  \| B^*(p_u - p_{uh}) \|_{X^*} \leq c_2^B \,
  \| p_u - p_{uh} \|_Y \leq c_2^B \inf\limits_{q_h \in Y_h}
  \| p_u - q_h \|_Y \, .
\end{equation}
In the same way we define $p_{fh} \in Y_h$ as unique solution of the
variational formulation
\[
  \langle A p_{fh} , q_h \rangle_H = \langle f , q_h \rangle_H \quad
  \mbox{for all} \; q_h \in Y_h ,
\]
in order to define $B^* p_{fh}$ as approximation of
$B^* p_f = B^* A^{-1} f$, i.e., $p_f = A^{-1} f$.
Hence, instead of \eqref{Galerkin VF Operator equation Su=g}
we now consider the perturbed variational formulation to find
$\widetilde{u}_H \in X_H$ such that
\begin{equation}\label{pert Galerkin VF Operator equation Su=g}
  \langle \widetilde{S} \widetilde{u}_H , v_H \rangle_H =
  \langle B^* p_{fh} , v_H \rangle_H \quad
  \mbox{for all} \; v_H \in X_H .
\end{equation}
To ensure unique solvability
of \eqref{pert Galerkin VF Operator equation Su=g},
we assume the discrete inf-sup stability condition
\begin{equation}\label{inf-sup Su=g}
  c_S \, \| v_H \|_X \leq
  \sup\limits_{0 \neq q_h \in Y_h}
  \frac{\langle B v_H , q_h \rangle_H}{\| q_h \|_Y}
  \quad \mbox{for all} \; v_H \in X_H .
\end{equation}

\begin{lemma}\label{lem:discrete-ellipticity}
  Assume the discrete inf-sup stability condition
  \eqref{inf-sup Su=g} to be satisfied. Then the approximate operator
  $\widetilde{S}$ is discrete elliptic in $X_H$, i.e., 
  \begin{equation}\label{discrete ellipticity S}
    \langle \widetilde{S} v_H , v_H \rangle_H \geq [c_S]^2 \,
    \| v_H \|_X^2 \quad \mbox{for all} \; v_H \in X_H .
  \end{equation}
\end{lemma}

\proof{We first note that $\widetilde{S} v_H = B^* p_{v_Hh}$ where
  $p_{v_Hh} \in Y_h$ solves
  \[
    \langle A p_{v_Hh} , q_h \rangle_H = \langle B v_H , q_h \rangle_H \quad
    \mbox{for all} \; q_h \in Y_h .
  \]
  Hence we obtain
  \[
    \langle \widetilde{S} v_H , v_H \rangle_H =
    \langle B^* p_{v_Hh} , v_H \rangle_H =
    \langle B v_H , p_{v_Hh} \rangle_H =
    \langle A p_{v_Hh} , p_{v_Hh} \rangle_H =
    \| p_{v_Hh} \|^2_Y \, .
  \]
  On the other hand, the discrete inf-sup stability condition
  \eqref{inf-sup Su=g} implies
  \[
    c_S \, \| v_H \|_X \leq
    \sup\limits_{0 \neq q_h \in Y_h}
    \frac{\langle B v_H , q_h \rangle_H}{\| q_h \|_Y} =
    \sup\limits_{0 \neq q_h \in Y_h}
    \frac{\langle A p_{v_Hh} , q_h \rangle_H}{\| q_h \|_Y} \leq
    \| p_{v_Hh} \|_Y ,
  \]
  and hence,
  \[
    c_S^2 \, \| v_H \|^2_X \, \leq \, \| p_{v_Hh} \|^2_Y \, = \,
    \langle \widetilde{S} v_H , v_H \rangle_H 
  \]
  follows.}

\noindent
The discrete ellipticity estimate \eqref{discrete ellipticity S} ensures
not only unique solvability of the perturbed variational formulation
\eqref{pert Galerkin VF Operator equation Su=g}, using the Strang lemma,
see, e.g., \cite{BrennerScott:2008,Steinbach:2008},
we can prove the following error estimate.

\begin{lemma}\label{Lemma Error abstract}
  For the unique solution $\widetilde{u}_H \in X_H$ of the perturbed
  variational formulation \eqref{pert Galerkin VF Operator equation Su=g}
  there holds the error estimate
  \begin{eqnarray}\label{Abstract error estimate}
  \| u -\widetilde{u}_H \|_X
  & \leq & \left( 1 + 2 \, \frac{[c_2^B]^2}{[c_S]^2} \right)
           \frac{c_2^B}{c_1^B} 
           \inf\limits_{v_H \in X_H} \| u - v_H \|_X \\
  && \hspace*{1cm} + \frac{c_2^B}{[c_S]^2} \Big[
     \inf\limits_{q_h \in Y_h} \| p_u - q_h \|_Y +
     \inf\limits_{q_h \in Y_h} \| p_f - q_h \|_Y \Big] \nonumber \, .
\end{eqnarray}
\end{lemma}

\proof{
  When considering the difference of the variational formulations
  \eqref{Galerkin VF Operator equation Su=g} and
  \eqref{pert Galerkin VF Operator equation Su=g}, this gives
  \[
    \langle S u_H - \widetilde{S} \widetilde{u}_H, v_H \rangle_H
    =
    \langle B^* (p_f-p_{fh}) , v_H \rangle_H \quad
    \mbox{for all} \; v_H \in X_H .
  \]
  From the discrete ellipticity
  \eqref{discrete ellipticity S} we then conclude
  \begin{eqnarray*}
    [c_S]^2 \, \| u_H - \widetilde{u}_H \|_X^2
    & \leq & \langle \widetilde{S} (u_H - \widetilde{u}_H) ,
             u_H - \widetilde{u}_H \rangle_H \\
    & = & \langle (\widetilde{S} - S)u_H , u_H - \widetilde{u}_H \rangle_H
          +
          \langle B^* (p_f - p_{fh}) , u_H - \widetilde{u}_H \rangle_H \\
    & \leq & \| (\widetilde{S}-S) u_H \|_{X^*}
             \| u_H - \widetilde{u}_H \|_X + c_2^B \,
             \| p_f - p_{fh} \|_Y \| u_H - \widetilde{u}_H \|_X,
  \end{eqnarray*}
  i.e.,
  \[
    [c_S]^2 \, \| u_H - \widetilde{u}_H \|_X
    \leq \| (\widetilde{S}-S) u_H \|_{X^*} +
    c_2^B \, \| p_f - p_{fh} \|_Y .
  \]
  We further have, using Lemma \ref{Lemma S} and \eqref{Stilde bounded},
  \begin{eqnarray*}
    \| (\widetilde{S}-S) u_H \|_{X^*}
    & \leq & \| (\widetilde{S}-S) u \|_{X^*} +
             \| (\widetilde{S}-S) (u-u_H) \|_{X^*} \\
    & \leq & \| (\widetilde{S}-S)u \|_{X^*} + 2 \, [c_2^B]^2 \,
             \| u - u_H \|_X ,
  \end{eqnarray*}
  and by using the triangle inequality,
  \[
    \| u -\widetilde{u}_H \|_X
    \leq \left( 1 + 2 \, \frac{[c_2^B]^2}{[c_S]^2} \right)
    \| u - u_H \|_X +
    \frac{1}{[c_S]^2} \Big[
    \| (\widetilde{S}-S) u \|_{X^*} +
    c_2^B \, \| p_f - p_{fh} \|_Y 
    \Big] .
  \]
  The assertion now follows from the norm equivalence inequalites
  \eqref{norm equivalence inequalities}, Cea's lemma
  \eqref{Cea Operator equation Su=g}, and the error estimate
  \eqref{Error Stilde}.}

\noindent
The perturbed variational formulation
\eqref{pert Galerkin VF Operator equation Su=g}
is, using $p_h:= p_{fh} - p_{\widetilde{u}_Hh} \in Y_h$, equivalent
to the coupled variational formulation to
find $(\widetilde{u}_H,p_h) \in X_H \times Y_h$ such that
\begin{equation}\label{Galerkin System abstract}
  \langle A p_h , q_h \rangle_H +
  \langle B \widetilde{u}_H , q_h \rangle_H
  = \langle f , q_h \rangle_H, \quad
  \langle p_h , B v_H \rangle_H = 0
\end{equation}
is satisfied for all $(v_H,q_h) \in X_H \times Y_h$. This is equivalent
to the linear system of algebraic equations,
\begin{equation}\label{eq:Linear system of saddle point formulation}
  \left( \begin{array}{cc}
           A_h & B_h \\
           B_h^\top &
         \end{array} \right)
       \left( \begin{array}{c}
                \underline{p} \\
                \underline{u}
              \end{array} \right)
            =
            \left( \begin{array}{c}
                     \underline{f} \\
                     \underline{0}
                     \end{array} \right),
\end{equation}
where, for $i,j=1,\ldots,M_Y$ and $k=1,\ldots,M_X$,
\[
  A_h[j,i] = \langle A \psi_i , \psi_j \rangle_H, \quad
  B_h[j,k] = \langle B \varphi_k , \psi_j \rangle_H, \quad
  f_j = \langle f , \psi_j \rangle_H .
\]
Since $A_h$ is symmetric and positive definite, and hence invertible,
we conclude the Schur complement system
\begin{equation}\label{LGS Schur}
B_h^\top A_h^{-1} B_h \underline{u} = B_h^\top A_h^{-1} \underline{f}
\end{equation}
which is the matrix representation of the perturbed variational
formulation \eqref{pert Galerkin VF Operator equation Su=g}.
Note that the Schur complement matrix $S_h := B_h^\top A_h^{-1} B_h$
is symmetric and positive definite.

Moreover, we observe that \eqref{Galerkin System abstract} is the
Galerkin formulation of the coupled variational formulation to
find $(u,p) \in X \times Y$ such that
\begin{equation}\label{VF Coupled system Su=g}
  \langle A p , q \rangle_H + \langle B u , q \rangle_H =
  \langle f , q \rangle_H, \quad
  \langle p , B v \rangle_H = 0 
\end{equation}
is satisfied for all $(v,q) \in X \times Y$, i.e., of the coupled
operator equation
\begin{equation}\label{Coupled system Su=g}
  A p + Bu = f , \quad B^* p = 0 \, .
\end{equation}
Due to $p := A^{-1}(f-Bu)$, this is \eqref{Operator equation Su=g}, where
by construction we have $p \equiv 0$. If the Galerkin matrix
$B_h$ is invertible, so is $B_h^\top$, and hence $p_h \equiv 0$ follows
in this particular case. But in general, the discrete inf-sup stability
condition \eqref{inf-sup Su=g} involves spaces such that
$\text{dim}(X_H)\neq \text{dim}(Y_h)$. Thus, $B_h$ is a not a square matrix
and hence not invertible, so that $p_h \in Y_h$ is not zero,
and we can use $p_h$ to define
an a posteriori error indicator for $\| u - \widetilde{u}_H \|_X$.

\begin{lemma}
  Let $(\widetilde{u}_H,p_h) \in X_H \times Y_h$ be the unique solution
  of \eqref{Galerkin System abstract}. Then,
  \begin{equation}\label{Estimate ph u-uH}
    \| p_h \|_Y \leq \| u - \widetilde{u}_H \|_S \leq c_2^B \,
    \| u - \widetilde{u}_H \|_X \, . 
  \end{equation}
\end{lemma}

\proof{When subtracting the Galerkin formulation
\eqref{Galerkin System abstract} from 
\eqref{VF Coupled system Su=g} for 
$q=q_h \in Y_h \subset Y$ and $v=v_H \in X_H \subset X$,
this gives the Galerkin orthogonalities
\begin{equation}\label{Galerkin orthogonality Coupled system Su=g}
  \langle A (p-p_h) , q_h \rangle_H +
  \langle B (u-\widetilde{u}_H) , q_h \rangle_H = 0, \quad
  \langle p - p_h , B v_H \rangle_H = 0
\end{equation}
for all $(v_H,q_h) \in X_H \times Y_h$.
In particular for $ p \equiv 0$ this gives
\[
  \langle A p_h, q_h \rangle_H = \langle B (u-\widetilde{u}_H), q	_h \rangle_H ,
  \quad
  \langle B v_H , p_h \rangle_H = 0 \quad
  \mbox{for all} \; (v_H,q_h) \in X_H \times Y_h .
\]
Hence, when choosing $q_h = p_h \in Y_h$, we further conclude
\[
  \| p_h \|_Y^2 = \langle A p_h , p_h \rangle_H =
  \langle B (u-\widetilde{u}_H) , p_h \rangle_H \leq 
  \| B (u - \widetilde{u}_H) \|_{Y^*} \| p_h \|_Y,
\]
i.e.,
\[
  \| p_h \|_Y \leq \| B(u - \widetilde{u}_H) \|_{Y^*} =
  \| u - \widetilde{u}_H \|_S ,
\]
and the assertion finally follows from \eqref{norm equivalence inequalities}.}

\noindent
While the upper estimate \eqref{Estimate ph u-uH} shows the efficiency of the
error estimator $\| p_h \|_Y$, reliability is more
involved. For this we introduce an ansatz space $X_{\overline{H}}\subset X$
such that $X_H \subset X_{\overline{H}}$ is satisfied. As in
\eqref{inf-sup Su=g} we assume the discrete inf-sup stability condition
\begin{equation}\label{inf-sup Xbar}
  \overline{c}_S \, \| v_{\overline{H}} \|_X \leq
  \sup\limits_{0 \neq q_h \in Y_h}
  \frac{\langle B v_{\overline{H}} , q_h \rangle_H}{\| q_h \|_Y}
  \quad \mbox{for all} \; v_{\overline{H}} \in X_{\overline{H}} .
\end{equation}
Due to $X_H \subset X_{\overline{H}}$ we have that
\eqref{inf-sup Su=g} is a direct consequence of \eqref{inf-sup Xbar}.
Using \eqref{inf-sup Xbar} we can determine
$(\widetilde{u}_{\overline{H}},\overline{p}_h) \in X_{\overline{H}} \times Y_h$
as unique solution of the variational formulation such that
\begin{equation}\label{VF Xbar}
  \langle A \overline{p}_h , q_h \rangle_H +
  \langle B \widetilde{u}_{\overline{H}} , q_h \rangle_H =
  \langle f , q_h \rangle_H, \quad
  \langle \overline{p}_h , B v_{\overline{H}} \rangle_H = 0
\end{equation}
is satisfied for all $(v_{\overline{H}},q_h) \in X_{\overline{H}} \times Y_h$.

\begin{lemma}
  Let $(\widetilde{u}_H,p_h) \in X_H \times Y_h$ and
  $(\widetilde{u}_{\overline{H}},\overline{p}_h) \in X_{\overline{H}} \times Y_h$
  be the unique solutions of the Galerkin variational formulations
  \eqref{Galerkin System abstract} and \eqref{VF Xbar}, respectively.
  Assume the saturation assumption
  \begin{equation}\label{saturation abstract}
    \| u - \widetilde{u}_{\overline{H}} \|_X \leq \eta \,
    \| u - \widetilde{u}_H \|_X \quad \mbox{for some} \; \eta \in (0,1).
  \end{equation}
  Then the error estimator $\| p_h \|_Y$ is reliable, satisfying
  \begin{equation}\label{reliable abstract}
    \| u - \widetilde{u}_H \|_X \leq \frac{1}{1-\eta} \,
    \frac{c_2^B}{\overline{c}_S^2} \,
    \| p_h \|_Y \, .
  \end{equation}
\end{lemma}

\proof{Subtracting the Galerkin variational formulation
  \eqref{VF Xbar} from \eqref{Galerkin System abstract}, this
  gives the Galerkin orthogonality
  \[
    \langle B (\widetilde{u}_{\overline{H}} - \widetilde{u}_H), q_h \rangle_H
    = \langle A (p_h - \overline{p}_h) , q_h \rangle_H \quad
    \mbox{for all} \; q_h \in Y_h .
  \]
  From the discrete inf-sup stability condition \eqref{inf-sup Xbar} we
  then conclude, recall
  $\widetilde{u}_H - \widetilde{u}_{\overline{H}} \in X_{\overline{H}}$,
  \[
    \overline{c}_S \, \| \widetilde{u}_H - \widetilde{u}_{\overline{H}} \|_X
    \leq \sup\limits_{0 \neq q_h \in Y_h}
    \frac{\langle B (\widetilde{u}_{\overline{H}} - \widetilde{u}_{H}),q_h
      \rangle_H}{\| q_h \|_Y}
    =
    \sup\limits_{0 \neq q_h \in Y_h}
    \frac{\langle A (p_h - \overline{p}_h) , q_h \rangle_H}{\| q_h \|_Y}
    \leq \| p_h - \overline{p}_h \|_Y \,. 
  \]
  On the other hand, using the above Galerkin orthogonality and the
  second equation in \eqref{VF Xbar} for
  $v_{\overline{H}} = \widetilde{u}_{\overline{H}} - \widetilde{u}_H$, this
  gives
  \begin{eqnarray*}
    \| p_h - \overline{p}_h \|_Y^2
    & = &
    \langle A (p_h - \overline{p}_h) , p_h - \overline{p}_h \rangle_H =
    \langle B(\widetilde{u}_{\overline{H}} - \widetilde{u}_H),
          p_h - \overline{p}_h \rangle_H \\
    & = & \langle B(\widetilde{u}_{\overline{H}} - \widetilde{u}_H),
          p_h \rangle_H  \leq c_2^B \,
          \| \widetilde{u}_{\overline{H}} - \widetilde{u}_H \|_X
          \| p_h \|_Y.
  \end{eqnarray*}
	\noindent
  Hence we obtain
  \[
    \| \widetilde{u}_{\overline{H}} - \widetilde{u}_H \|_X^2 \leq
    \frac{1}{\overline{c}_S^2} \,
    \| p_h - \overline{p}_h \|^2_Y \leq \frac{c_2^B}{\overline{c}_S^2} \,
    \| \widetilde{u}_{\overline{H}} - \widetilde{u}_H \|_X
    \| p_h \|_Y ,
  \]
  i.e.,
  \[
    \| \widetilde{u}_{\overline{H}} - \widetilde{u}_H \|_X \leq
    \frac{c_2^B}{\overline{c}_S^2} \,
    \| p_h \|_Y \, .
  \]
  Using the triangle inequality and the saturation assumption
  \eqref{saturation abstract} we finally have
  \[
    \| u - \widetilde{u}_H \|_X \leq
    \| u - \widetilde{u}_{\overline{H}} \|_X +
    \| \widetilde{u}_{\overline{H}} - \widetilde{u}_H \|_X \leq
    \eta \, \| u - \widetilde{u}_H \|_X + \frac{c_2^B}{\overline{c}_S^2} \,
    \| p_h \|_Y ,
  \]
  from which the assertion follows.}

\noindent
It remains to define, for a given ansatz space $X_{\overline{H}}$, the test
space $Y_h$ such that the discrete inf-sup condition \eqref{inf-sup Xbar}
and therefore \eqref{inf-sup Su=g} are satisfied. This can be achieved by
assuming a sufficiently rich test space $Y_h$, as stated in the following
abstract approach.  But, as we will see later,
 this is not always required, since one may establish
\eqref{inf-sup Xbar} in a different way.

\begin{theorem}
  For a given finite element space $X_{\overline{H}} \subset X$ let
  $Y_h \subset Y$ such that
  \begin{equation}\label{Assumption Yh}
    \sup\limits_{v_{\overline{H}} \in X_{\overline{H}}} \inf\limits_{q_h \in Y_h}
    \| p_{v_{\overline{H}}} - q_h \|_Y \leq \delta \,
    \| p_{v_{\overline{H}}} \|_Y = \delta \,
    \| v_{\overline{H}} \|_S
  \end{equation}
  is satisfied for some $\delta \in (0,1)$. Then there holds the
  discrete inf-sup stability condition \eqref{inf-sup Xbar}, i.e.,
  \begin{equation}\label{inf sup delta}
    c_1^B (1-\delta) \, \| v_{\overline{H}} \|_X \leq
    \sup\limits_{q_h \in Y_h} \frac{\langle B v_{\overline{H}},q_h \rangle_H}
    {\| q_h \|_Y} \quad \mbox{for all} \;
    v_{\overline{H}} \in X_{\overline{H}} .
  \end{equation}
  
\end{theorem}

\proof{For an arbitrary but fixed $v_{\overline{H}} \in X_{\overline{H}}$ we
  define $p_{v_{\overline{H}}} = A^{-1} B v_{\overline{H}} \in Y$ and as in
  the proof of Lemma \ref{Lemma S} we conclude
  \[
    \| p_{v_{\overline{H}}} \|^2_Y =
    \| v_{\overline{H}} \|_S^2 =
    \langle B v_{\overline{H}} , p_{v_{\overline{H}}} \rangle_H \, .
  \]
  In addition we define $p_{v_{\overline{H}}h} \in Y_h$ as unique
  solution of the Galerkin variational formulation
  \[
    \langle A p_{v_{\overline{H}}h},q_h \rangle_H =
    \langle B v_{\overline{H}} , q_h \rangle_H =
    \langle A p_{v_{\overline{H}}},q_h \rangle_H
    \quad \mbox{for all} \; q_h \in Y_h ,
  \]
  satisfying the bound
  \[
    \| p_{v_{\overline{H}}h} \|_Y \leq \| p_{v_{\overline{H}}} \|_Y,
  \]
  and Cea's lemma,
  \[
    \| p_{v_{\overline{H}}} - p_{v_{\overline{H}}h} \|_Y \leq
    \inf\limits_{q_h \in Y_h} \| p_{v_{\overline{H}}} - q_h \|_Y \, .
  \]
  From \eqref{Assumption Yh} we then obtain
  \[
    \| p_{v_{\overline{H}}} - p_{v_{\overline{H}}h} \|_Y \leq \delta \,
    \| p_{v_{\overline{H}}} \|_Y \, .
  \]
  Hence we can write
  \begin{eqnarray*}
    \langle B v_{\overline{H}} , p_{v_{\overline{H}}h} \rangle_H
    & = & \langle A p_{v_{\overline{H}}} , p_{v_{\overline{H}}h} \rangle_H \, = \,
          \langle A p_{v_{\overline{H}}} , p_{v_{\overline{H}}} \rangle_H -
          \langle A p_{v_{\overline{H}}} ,
          p_{v_{\overline{H}}} - p_{v_{\overline{H}}h} \rangle_H \\
    & \geq & \| p_{v_{\overline{H}}} \|_Y^2 -
             \| p_{v_{\overline{H}}} \|_Y 
             \| p_{v_{\overline{H}}} - p_{v_{\overline{H}}h} \|_Y \, \geq \,
             (1-\delta) \, \| p_{v_{\overline{H}}} \|_Y^2 \\
    & \geq & (1-\delta) \, \| p_{v_{\overline{H}}} \|_Y
             \| p_{v_{\overline{H}}h} \|_Y, 
  \end{eqnarray*}
  i.e.,
  \[
    (1-\delta) \, \| v_{\overline{H}} \|_S =
    (1-\delta) \, \| p_{v_{\overline{H}}} \|_Y \leq
    \frac{\langle B v_{\overline{H}} , p_{v_{\overline{H}}h} \rangle_H}
    {\| p_{v_{\overline{H}}h} \|_Y},
  \]
  implying the inf-sup stability condition \eqref{inf sup delta}.}

\noindent
In some applications, e.g., when considering space-time finite element
methods for the heat equation as in \cite{Steinbach:2015},
the discrete inf-sup condition can be established when using a 
(discretization dependent) norm on the ansatz space $X$. In particular,
let $\norm{\cdot}_{X,h} : X \to \mathbb{R}$ define a norm on $X_H$, satisfying 
$\| v \|_{X,h} \leq \| v \|_X$ for all $v \in X$ and assume that 
\begin{equation}\label{eq:discrete-inf-sup-XH}
  \widetilde{c}_S \, \| v_H \|_{X,h} \leq
  \sup_{0 \neq q_h \in Y_h} \frac{\langle Bv_H,q_h\rangle_H}
  {\| q_h \|_Y} \quad \mbox{for all} \; v_H \in X_H.
\end{equation} 
Then the following stability and error estimates hold true. 

\begin{lemma}\label{lem:error-estimate-XH}
  Assume the discrete inf-sup stability condition
  \eqref{eq:discrete-inf-sup-XH}. Then the approximate operator
  $\widetilde S$ is discrete elliptic, satisfying
  \begin{equation*}
    \langle \widetilde{S} v_H,v_H \rangle_H \geq
    [\widetilde{c}_S]^2 \, \| v_H \|_{X,h}^2 \quad
    \mbox{for all} \; v_H \in X_H,
  \end{equation*}
  and the perturbed variational formulation
  \eqref{pert Galerkin VF Operator equation Su=g} admits a unique solution.
  Furthermore, there holds the error estimate
  \begin{equation}\label{eq:error-estimate-XH}
    \| u-\widetilde{u}_H \|_{X,h} \leq
    \Big( 1+\frac{2c_2^B}{\widetilde{c}_S} \Big)
    \inf_{v_H \in X_H} \| u-v_H \|_X . 
  \end{equation}
\end{lemma}

\proof{The proof of the discrete ellipticity estimate follows the lines of
  the proof of Lemma \ref{lem:discrete-ellipticity}, replacing the discrete
  inf-sup stability condition \eqref{inf-sup Su=g} by
  \eqref{eq:discrete-inf-sup-XH}. This already ensures unique solvability
  of \eqref{pert Galerkin VF Operator equation Su=g}, as
  $\| \cdot \|_{X,h}$ defines a norm on $X_H \subset X$. To derive the error
  estimate, let $v_H \in X_H$ be arbitrary but fixed. First note that we have 
  \begin{equation*}
    \| u-\widetilde{u}_H \|_{X,h} \leq
    \| u-v_H \|_{X,h} + \| v_H-\widetilde{u}_H \|_{X,h} \leq
    \| u-v_H \|_X + \| v_H - \widetilde{u}_H \|_{X,h} \, .
  \end{equation*}
  When using the discrete inf-sup stability condition
  \eqref{eq:discrete-inf-sup-XH}, and \eqref{Galerkin System abstract},
  for the second term we further have
  \begin{eqnarray*}
    \widetilde{c}_S \, \| v_H - \widetilde{u}_H \|_{X,h}
    & \leq & \sup_{0\neq q_h \in Y_h}
             \frac{\langle B(v_H-\widetilde{u}_H),q_h\rangle_H}
             {\| q_h \|_Y} \, = \,
             \sup_{0\neq q_h\in Y_h} \frac{\langle Bv_H-(f-Ap_h),q_h\rangle_H}
             {\| q_h \|_Y} \\
    & = & \sup_{0\neq q_h\in Y_h}
          \frac{\langle B(v_H-u)+Ap_h,q_h\rangle_H}{\| q_h \|_Y} \, \leq \,
          c_2^B \, \| u-v_H \|_X + \| p_h \|_Y .
  \end{eqnarray*}
  Using $\langle p_h,Bv_H\rangle_H=0$ for all $v_H\in X_H$, see
  \eqref{Galerkin System abstract}, we can further estimate
  \[
    \| p_h \|_Y^2 = \langle Ap_h,p_h \rangle_H =
  \langle B(u-\widetilde{u}_H),p_h \rangle_H =
  \langle B(u-v_H),p_h\rangle_H \leq
  c_2^B \, \| u-v_H \|_X \| p_h \|_Y,  
\]
i.e., $\| p_h \|_Y \leq c_2^B \, \| u-v_H \|_X$. Since $v_H\in X_H$ was
arbitrary, this concludes the proof.}

\noindent
In what follows we will discuss several applications of this abstract
setting. Although our main interest is in the discretization of time dependent
partial differential equations such as the heat and the wave equation,
we will first consider an elliptic problem in order to present
the main ideas of this approach for examples which are well known in
literature. Crucial is the choice of the finite element spaces
$X_{\overline{H}}$ and $Y_h$ such that the discrete inf-sup condition
\eqref{inf-sup Xbar} is satisfied. As we will see, and depending on the
particular partial differential equation to be solved, we may even
consider $X_{\overline{H}}=Y_h$, or we may consider stable pairs of finite
element functions which are defined with respect to the same finite
element mesh. In the most general case we may choose first a finite
element space $X_H$, then a space $X_{\overline{H}}$ to ensure the
saturation condition \eqref{saturation abstract}, and finally
$Y_h$ to satisfy the discrete inf-sup stability condition
\eqref{inf-sup Xbar}, e.g., we may define
the ansatz space $Y_h$ with respect to a sufficient small mesh size
$h < c_0 \overline{H}$ in order to meet \eqref{Assumption Yh}.

\section{Elliptic Dirichlet boundary value problem}
As first example we consider the Dirichlet boundary value problem for
the Poisson equation,
\begin{equation}\label{DBVP Poisson}
  - \Delta u(x) = f(x) \quad \mbox{for} \; x \in \Omega, \quad
  u(x) = 0 \quad \mbox{for} \; x \in \Gamma,
\end{equation}
where $\Omega \subset {\mathbb{R}}^n$, $n=2,3$, is a bounded Lipschitz
domain with boundary $\Gamma=\partial \Omega$. With respect to the
abstract setting we have
\[
  X = Y = H^1_0(\Omega), \quad
  A = B = - \Delta : H^1_0(\Omega) \to H^{-1}(\Omega), \quad
  \| v \|_X = \| \nabla v \|_{L^2(\Omega)} .
\]
We obviously have \eqref{Assumptions B} with $c_1^B=c_2^B=1$.
In this case, the abstract variational formulation
\eqref{VF Coupled system Su=g} reads to find
$(u,p) \in H^1_0(\Omega) \times H^1_0(\Omega)$ such that
\begin{equation}\label{DBVP Poisson VF}
  \int_\Omega \nabla p(x) \cdot \nabla q(x) \, dx +
  \int_\Omega \nabla u(x) \cdot \nabla q(x) \, dx =
  \langle f , q \rangle_\Omega, \quad
  \int_\Omega \nabla v(x) \cdot \nabla p(x) \, dx = 0 \quad
\end{equation}
is satisfied for all $(v,q) \in H^1_0(\Omega) \times H^1_0(\Omega)$. Let
$X_h :=S_h^1(\Omega)\cap H_0^1(\Omega)=\text{span}\{\varphi_k\}_{k=1}^{M_X}$,
where $S_h^1(\Omega)=\text{span}\{\varphi_k\}_{k=1}^{\widetilde{M}_X}$
denotes the standard finite element ansatz space
of piecewise linear and continuous basis functions, which are defined with
respect to some admissible locally quasi-uniform decomposition of
$\Omega$ into shape regular simplicial finite elements $\tau_\ell^h$ of
local mesh size $h_\ell$. By definition we have, for $ 0 \neq u_h \in X_h$,
\[
  \| \nabla u_h \|_{L^2(\Omega)} =
  \frac{\langle \nabla u_h , \nabla u_h \rangle_{L^2(\Omega)}}
  {\| \nabla u_h \|_{L^2(\Omega)}} \leq
  \sup\limits_{0 \neq q_h \in X_h}
  \frac{\langle \nabla u_h , \nabla q_h \rangle_{L^2(\Omega)}}
  {\| \nabla q_h \|_{L^2(\Omega)}} ,
\]
i.e., \eqref{inf-sup Xbar} for $Y_h = X_{\overline{H}}=X_h$.
In the case $X_H=X_{\overline{H}}= X_h$ we finally obtain, due to $p_h=0$, the
standard finite element formulation for the Poisson equation, i.e.,
\begin{equation}\label{standard FEM}
  \int_\Omega \nabla u_h(x) \cdot \nabla q_h(x) \, dx =
  \langle f , q_h \rangle_\Omega \quad \mbox{for all} \; q_h \in X_h .
\end{equation}
In order to use $\| \nabla p_h \|_{L^2(\Omega)}$ as a reliable and
efficient error estimator we need to introduce a finite
element space $X_H$ which is defined with respect to a coarser mesh.
In fact, one may first define the ansatz space $X_H$, and afterwards
$Y_h$ by applying some additional refinements. Hence we consider
the mixed variational formulation to find
$(\widetilde{u}_H , p_h) \in X_H \times Y_h$ such that
\begin{equation}\label{DBVP Poisson FEM}
  \int_\Omega \nabla p_h(x) \cdot \nabla q_h(x) \, dx +
  \int_\Omega \nabla \widetilde{u}_H(x) \cdot \nabla q_h(x) \, dx =
  \langle f , q_h \rangle_\Omega , \quad
  \int_\Omega \nabla v_H(x) \cdot \nabla p_h(x) \, dx = 0  
\end{equation}
is satisfied for all $(v_H,q_h) \in X_H \times Y_h$. According
to Lemma \ref{Lemma Error abstract} we conclude the error estimate
\[
  \| \nabla (u-\widetilde{u}_H) \|_{L^2(\Omega)} \leq
  c_1 \inf\limits_{v_H \in X_H} \| \nabla (u-v_H) \|_{L^2(\Omega)} +
  c_2 \inf\limits_{q_h \in Y_h} \| \nabla (u-q_h) \|_{L^2(\Omega)} .
\]
Hence, assuming $u \in H^1_0(\Omega) \cap H^s(\Omega)$ for some
$ s \in [1,2]$ we finally obtain the error estimate
\begin{equation}\label{FEM error Laplace}
  \| \nabla (u-\widetilde{u}_H) \|_{L^2(\Omega)} \leq
  c \, H^{s-1} \, |u|_{H^s(\Omega)} .
\end{equation}
Note that the saturation assumption \eqref{saturation abstract} now reads,
for $u_h \in X_h = Y_h$,
\[
  \| \nabla (u - u_h) \|_{L^2(\Omega)} \leq \eta \,
  \| \nabla (u - \widetilde{u}_H) \|_{L^2(\Omega)} \quad
  \mbox{for some} \; \eta \in (0,1).
\]
For the particular choice $H=2h$, i.e., one additional refinement
to define $Y_h$ from $X_H$, this is obviously related
to the $h-\frac{h}{2}$ error estimator, e.g.,
\cite{Praetorius:2014,
  2010_FerrazOrtnerPraetorius_ConvergenceofSimpleAdaptiveGalerkin}.
Although this error estimator is well established in literature,
we just present one example solving the mixed Galerkin finite
element scheme \eqref{DBVP Poisson FEM}. Note that in contrast
to the standard $h-\frac{h}{2}$ error estimator we compute
$(p_h,\widetilde{u}_H)$ at once solving \eqref{DBVP Poisson FEM}.

In what follows we will use the error estimator
\[
  \eta_h^2 := \| \nabla p_h \|_{L^2(\Omega)}^2 = \int_\Omega
  |\nabla p_h(x)|^2 \, dx =
  \sum\limits_{\ell=1}^{N_h} \int_{\tau_\ell^h} |\nabla p_h(x)|^2 \, dx =
  \sum\limits_{\ell=1}^{N_h} \eta_{h,\ell}^2
\]
with the local error indicators
\[
\eta_{h,\ell}^2 = \int_{\tau_\ell^h} |\nabla p_h(x)|^2 \, dx \, .
\]
In order to define local indicators for the error
$\| \nabla (u-\widetilde{u}_H) \|_{L^2(\Omega)}$ we sum up
all local contributions,
\[
\eta_{H,j}^2 = \sum\limits_{\tau_\ell^h \subset \tau_j^H} \eta_{h,\ell}^2 \, .
\]
For the adaptive refinement we use the D\"orfler marking 
\cite{Doerfler:1996} with the parameter $\theta=0.5$.
The marked elements are then refined using 
newest vertex bisection.
All computations were done using the finite element software Netgen/NGSolve
\cite{Schoeberl:2014}. The resulting linear systems were solved using
the sparse direct solver package Umfpack.

\begin{figure}[!htbp]
	\centering
	\begin{tikzpicture}
		\begin{axis}[
			xmode = log,
			ymode = log,
			xlabel= $\widetilde{M}_X$,
			ylabel= {$\norm{\nabla p_h}_{L^2(\Omega)}$, $\norm{\nabla (u-\widetilde{u}_H)}_{L^2(\Omega)}$},
			ymin = 0.001,
			ymax = 1,
			legend style={font=\tiny}, legend pos = north east]
			\addplot[mark = x,blue] table [col sep=&, y=L2_dxph, x=nv]{tables/SimData_LSLaplace_Lshape_adaptive.dat};
			\addlegendentry{$\norm{\nabla p_h}_{L^2(\Omega)}$}
			
			\addplot[mark = x,red] table [col sep=&, y=errH1_u, x=nv]{tables/SimData_LSLaplace_Lshape_adaptive.dat};
			\addlegendentry{$\norm{\nabla (u-\widetilde{u}_H)}_{L^2(\Omega)}$}
			
		\end{axis}
	\end{tikzpicture}
	\caption{Comparison of the error $\| \nabla (u-\widetilde{u}_H)
          \|_{L^2(\Omega)}$ and the error estimator
          $\| \nabla p_h\|_{L^2(\Omega)}$.}
	\label{abb:laplace_ex1_errorVSestimator}
\end{figure}
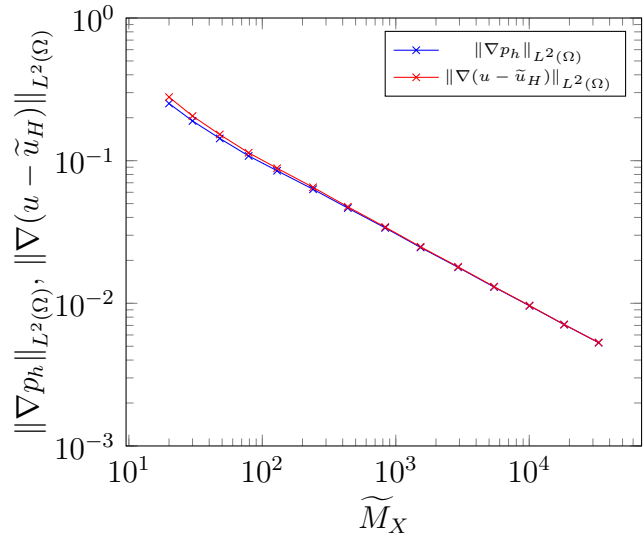

As example we consider the $L$ shaped domain 
\[
  \Omega = (-1,1)^2 \setminus
  \left(\left[0,1\right] \times \left[-1,0\right]\right)
  \subset \R^2,
\]
and we consider the Dirichlet boundary value problem
\eqref{DBVP Poisson} with the exact solution, given in polar coordinates,
\begin{equation}\label{singular example}
  u(r,\varphi) = r^{2/3} \, \sin \frac{2}{3}\varphi, \quad
  u \in H^s(\Omega), \; s < \frac{5}{3} .
\end{equation}
Due to the reduced regularity we expect a reduced order of 
convergence of $\mathcal{O}(H^{2/3})$ for a uniform refinement
when measuring the energy error $\| \nabla (u-\widetilde{u}_H) \|_{L^2(\Omega)}$,
and of $\mathcal{O}(H^{4/3})$ when considering
the $L^2$ error $ \| u - \widetilde{u}_H \|_{L^2(\Omega)}$. When using
the adaptive refinement strategy as described, we recover the
optimal convergence of ${\mathcal{O}}(H)$ in the energy norm and
of ${\mathcal{O}}(H^2)$ in $L^2(\Omega)$, respectively. The related
numerical results are given in Fig.~\ref{abb:laplace_ex1_errors}.
In Fig.~\ref{abb:laplace_ex1_errorVSestimator} we present a comparison
of the error $\norm{\nabla (u-\widetilde{u}_H)}_{L^2(\Omega)}$ and the error estimator
$\eta=\norm{\nabla p_h}_{L^2(\Omega)}$ which shows that the error estimator is effective.
			
\begin{figure}[!htbp]
	\centering
	\begin{tikzpicture}[scale=1]
		\begin{axis}[
			xmode = log,
			ymode = log,
			xlabel= $\widetilde{M}_X$,
			ylabel= errors,
			ymin = 0.00001,
			ymax = 1,
			legend style={font=\tiny}, legend pos = outer north east]
			\addplot[mark = *,dashed,blue] table [col sep=&, y=errL2_u, x=nv]{tables/SimData_LSLaplace_Lshape_uniform.dat};
			\addlegendentry{$\norm{u-\widetilde{u}_H}_{L^2(\Omega)}$ uniform}
			
			\addplot[mark = *,dashed,red] table [col sep=&, y=errH1_u, x=nv]{tables/SimData_LSLaplace_Lshape_uniform.dat};
			\addlegendentry{$\norm{\nabla (u-\widetilde{u}_H)}_{L^2(\Omega)}$ uniform}
			
			\addplot[mark = square*,blue] table [col sep=&, y=errL2_u, x=nv]{tables/SimData_LSLaplace_Lshape_adaptive.dat};
			\addlegendentry{$\norm{u-\widetilde{u}_H}_{L^2(\Omega)}$ adaptive}
			
			\addplot[mark = square*,red] table [col sep=&, y=errH1_u, x=nv]{tables/SimData_LSLaplace_Lshape_adaptive.dat};
			\addlegendentry{$\norm{\nabla (u-\widetilde{u}_H)}_{L^2(\Omega)}$ adaptive}
			\addplot[
				domain = 100:200000,
				samples = 10,
				dashed,
				thin,
				red,
			] {x^(-1/3)};
			\addlegendentry{$\widetilde{M}_X^{-1/3}\sim H^{2/3}$}
			\addplot[
				domain = 100:200000,
				samples = 10,
				dashed,
				thin,
				blue,
			] {0.15*x^(-2/3)};
			\addlegendentry{$\widetilde{M}_X^{-2/3}\sim H^{4/3}$}
			
			\addplot[
				domain = 100:10000,
				samples = 10,
				thin,
				red,
			] {0.7*x^(-1/2)};
			\addlegendentry{$\widetilde{M}_X^{-1/2}\sim H$}
			
			\addplot[
				domain = 100:10000,
				samples = 10,
				thin,
				blue,
			] {0.3*x^(-1)};
			\addlegendentry{$\widetilde{M}_X^{-1}\sim H^{2}$}
			
		\end{axis}
	\end{tikzpicture}
	\caption{Convergence behavior for the approximation of the
          singular function \eqref{singular example}.}
	\label{abb:laplace_ex1_errors}
\end{figure}
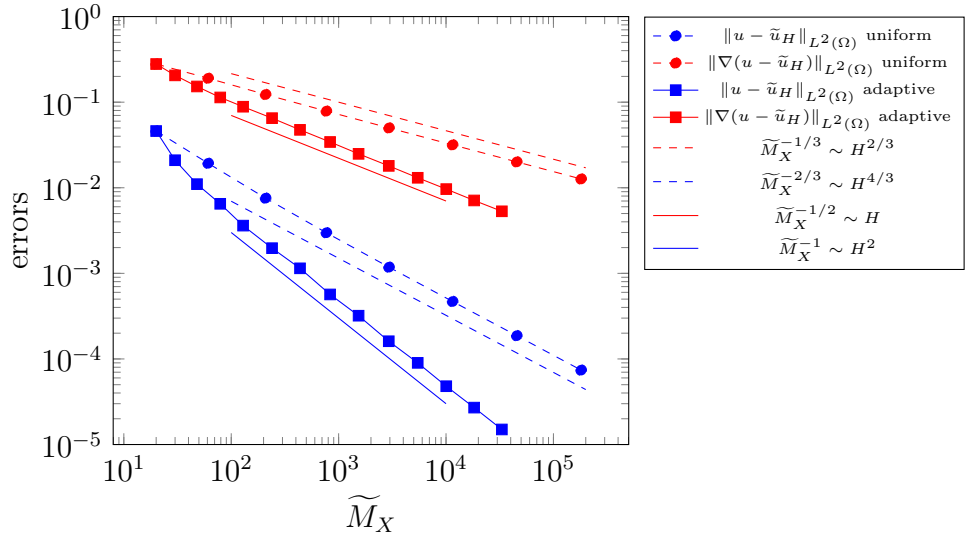	
      
\begin{remark}
  Instead of $X=Y=H^1_0(\Omega)$, $S = B^* A^{-1} B = - \Delta$,
  we may also consider the abstract setting with
  $X = H_\Delta(\Omega) := \{ v \in H^1_0(\Omega) :
  \Delta v \in L^2(\Omega)\}$, $Y = L^2(\Omega)$, 
  $B = - \Delta : H_\Delta(\Omega) \to Y$, and
  $A = I : L^2(\Omega) \to L^2(\Omega)$. Then we have that
  $S = B^* A^{-1} B = \Delta^2$ is the Bi-Laplacian. 
  For a conforming finite element
  discretization in this case one may use tensor product meshes
  and quadratic B splines to define $X_h$, and piecewise constant
  basis functions to define $Y_h$. For a related approach in the case
  of a distributed optimal control problem, see, e.g.,
  \cite{Brenner:2020, Brenner:2023}. Alternatively, and applying integration
  by parts twice, one may exchange the roles of $X=L^2(\Omega)$ and
  $Y = H_\Delta(\Omega)$, in order to consider, e.g., less regular
  Dirichlet boundary conditions $u=g \in L^2(\Gamma)$, see, e.g.,
  \cite{Apel:2017}. In addition, one can also apply the abstract theory
  as presented in this paper
  for the analysis of related boundary element methods
  \cite{Steinbach:LSBEM}.
\end{remark}

\section{Parabolic Dirichlet boundary value problem}
As an example for a parabolic partial differential equation we
consider the Dirichlet problem for the heat equation,
\begin{subequations}\label{eq:Parabolic DBVP}
  \begin{alignat}{2}
    \partial_t u(x,t) -\Delta_x u(x,t)
    & = f(x,t) & \quad & \text{for }(x,t) \in Q:=\Omega \times (0,T), \\
    u(x,t) & = 0 & \quad & \text{for }(x,t)\in \Sigma:=\Gamma \times (0,T),\\
    u(x,0) & = 0 & \quad & \text{for } x \in \Omega,
  \end{alignat}
\end{subequations}
where $\Omega\subset \R^n$, $n=1,2,3$ is a bounded Lipschitz domain
with boundary $\Gamma =\partial \Omega$, and $T>0$ is a given time horizon.
In view of the abstract setting we have the Bochner spaces
\[
  X:= L^2(0,T;H_0^1(\Omega))\cap H_{0,}^1(0,T;H^{-1}(\Omega)),~
  Y:= L^2(0,T;H_0^1(\Omega)),~
  H: = L^2(Q),
\]
where $H_{0,}^1(0,T;H^{-1}(\Omega)) := \{ u \in L^2(Q) :
\partial_tu \in L^2(0,T;H^{-1}(\Omega)), u(x,0)=0, x \in \Omega \}$.
The operators $A:Y \to Y^*$ and $B:X\to Y^*$ are defined in a variational
sense satisfying
\[
  \langle Ap,q \rangle_Q := \langle \nabla_x p , \nabla_x q \rangle_{L^2(Q)},
  \quad
  \langle Bu,q \rangle_Q := \langle \partial_t u,q \rangle_Q +
  \langle \nabla_x u , \nabla_x q \rangle_{L^2(Q)}
\]
for all $p,q \in Y$ and $u \in X$. The corresponding norms read
\[
  \| p \|_Y := \| \nabla_x p \|_{L^2(Q)}, \quad 
  \| u \|_X := \sqrt{ \| \partial_t u \|_{Y^*}^2 + \| \nabla_x u \|_Y^2},
  \quad
  \| \partial_t u \|_{Y^*} = \| \nabla_x w_u \|_{L^2(Q)},
\]
where $w_u \in Y$ is the unique solution of the variational problem
\begin{equation}\label{Def wu}
  \langle \nabla_x w_u , \nabla_x q \rangle_{L^2(Q)} =
  \langle \partial_t u , q \rangle_Q \quad \mbox{for all} \; q \in Y.
\end{equation}
The operator $A : Y \to Y^*$ is self-adjoint, bounded and elliptic.
Moreover, $B : X \to Y^*$ is bounded with $c_2^B=\sqrt{2}$, i.e.,
\[
  \| B v \|_{Y^*} \leq \sqrt{2} \, \| v \|_X \quad
  \mbox{for all} \; v \in X .
\]
In order to prove an inf-sup stability condition, for $u \in X$ we first define
$ w_u \in Y$ as in \eqref{Def wu}. For $q := u + w_u \in Y$ we then have
\begin{eqnarray*}
  \spf{Bu}{q}_Q
  & = & \spf{Bu}{u+w_u}_Q \\
  & = & \langle \partial_t u , u \rangle_Q +
        \langle \nabla_x u , \nabla_x u \rangle_{L^2(Q)} +
        \langle \partial_t u , w_u \rangle_Q +
        \langle \nabla_x u , \nabla_x w_u \rangle_{L^2(Q)} \\
  & = & \langle \nabla w_u , \nabla_x u \rangle_{L^2(\Omega)} +
        \langle \nabla_x u , \nabla_x u \rangle_{L^2(Q)} +
        \langle \nabla_x w_u , \nabla_x w_u \rangle_{L^2(Q)} +
        \langle \nabla_x u , \nabla_x w_u \rangle_{L^2(Q)} \\
  & = & \langle \nabla_x (u+w_u) , \nabla_x (u+w_u) \rangle_{L^2(\Omega)} \\
  & = & \| q \|_Y^2 \, .
\end{eqnarray*}
On the other hand,
\begin{eqnarray*}
  \| q \|^2_Y 
  & = & \| \nabla (u + w_u) \|_{L^2(Q)}^2 \\
  & = & \langle \nabla_x u , \nabla_x u \rangle_{L^2(\Omega)} +
        \langle \nabla_x w_u , \nabla_x w_u \rangle_{L^2(Q)} +
        2 \,\langle \nabla_x w_u , \nabla_x u \rangle_{L^2(Q)} \\
  & = & \langle \nabla_x u , \nabla_x u \rangle_{L^2(\Omega)} +
        \langle \nabla_x w_u , \nabla_x w_u \rangle_{L^2(Q)} +
        2 \, \langle \partial_t u , u \rangle_Q \\
  & \geq & \| \nabla_x u \|_Y^2 + \| \nabla_x w_u \|^2_Y \\
  & = & \| u \|_X^2 ,
\end{eqnarray*}
implies
\[
\spf{Bu}{u+w_u}_Q \geq \| u \|_X \| u + w_u \|_Y,
\]
and therefore the inf-sup stability condition
\[
  \| u \|_X \leq \sup\limits_{0 \neq q \in Y} \frac{\spf{Bu}{q}_Q}{\| q \|_Y}
  \quad \mbox{for all} \; u \in X
\]
follows with $c_1^B = 1$. Note that this estimate is improved than
originally derived in \cite{Steinbach:2015}.
Moreover, $B:X \to Y^*$ is surjective, see, e.g.,
\cite{LSTY:2021}. Thus, we can apply the abstract theory from
Section \ref{sec:abstract setting}.
The variational formulation \eqref{VF Coupled system Su=g} now reads to
find $(u,p)\in X \times Y$ such that 
\begin{eqnarray}
  \int_Q \nabla_x p(x,t) \cdot \nabla_x q(x,t) \, dx \, dt
  \hspace*{6cm} && \nonumber \\
  + \int_Q \Big[ \partial_t u(x,t) \, q(x,t) +
    \nabla_x u(x,t) \cdot \nabla_x q(x,t) \Big] dx \, dt
  & = & \langle f,q \rangle_Q, \label{eq:DBVP Parabolic VF} \\
  \int_Q \Big[ \partial_t v(x,t) \, p(x,t) + \nabla_x v(x,t) \cdot
    \nabla_x p(x,t) \Big] dx \, dt
    & = & 0 \nonumber
\end{eqnarray}
is satisfied for all $(v,q) \in X \times Y$.

For the discretization of \eqref{eq:DBVP Parabolic VF} we first consider
an ansatz space $X_H= \mbox{span} \{ \varphi_k \}_{k=1}^{M_X} \subset X$
of piecewise linear and continuous basis functions $\varphi_k$ which are
defined with respect to some admissible and locally quasi-uniform
decomposition of the space-time domain $Q$ into shape regular simplicial
finite elements $q_\ell^H$ of mesh size $H_\ell$. In addition we define
the ansatz space $Y_h = \mbox{span} \{ \psi_i \}_{i=1}^{M_Y} \subset Y$
of piecewise linear and continuous basis functions $\psi_i$ which are
defined with respect to some refined decomposition of the space-time
domain $Q$ into finite elements $q_i^h$ of local mesh size $h_i$.
From a practical point of view, we may use one additional refinement
of the mesh which was used to define $X_H$, i.e., for $q_i^h \subset q_\ell^H$
we have $h_i = \frac{1}{2}H_\ell$. As an alternative, we may also use
second order form functions on $q_\ell^H$ to define $Y_h$. In both cases
we have the inclusion $X_H \subset Y_h$ which enables us to prove
a discrete inf-sup stability condition.
For any $u \in X$ we define $w_{uh} \in Y_h$ as the unique solution
of the Galerkin variational formulation
\[
  \langle \nabla_x w_{uh} , \nabla_x q_h \rangle_{L^2(Q)} =
  \langle \partial_t u , q_h \rangle_Q \quad \mbox{for all} \; q_h \in Y_h.
\]
With this we define the discrete norm
\[
  \| u \|_{X,h} := \sqrt{ \| u \|_Y^2 + \| w_{uh} \|_Y^2} \leq \,
  \| u \|_X ,
\]
and as in the continuous case, instead of $w_u\in Y$ we now use
$w_{u_Hh} \in Y_h$, we can prove the discrete inf-sup stability condition
\[
  \| u_H \|_{X,h} \leq \sup\limits_{0 \neq q_h \in Y_h}
  \frac{\spf{Bu_H}{q_h}_Q}{\|q_h\|_Y} \quad \mbox{for all} \;
  u_H \in X_H .
\]
Hence we conclude unique solvability of the mixed space-time
finite element variational formulation to find
$(\widetilde{u}_H,p_h) \in X_H \times Y_h$ such that
\begin{eqnarray}
  \int_Q \nabla_x p_h(x,t) \cdot \nabla_x q_h(x,t) \, dx \, dt
  \hspace*{6cm} \nonumber
  && \\
  + \int_Q \Big [\partial_t \widetilde{u}_H(x,t) \, q_h(x,t)  +
  \nabla_x \widetilde{u}_H(x,t) \cdot \nabla_x q_h(x,t) \Big] dx \, dt
  &=& \langle f , q_h \rangle_Q, \label{eq:DBVP Parabolic FEM} \\
  \int_Q \Big[ \partial_t v_H(x,t) \, p_h(x,t) +
  \nabla_x v_H(x,t) \cdot \nabla_x p_h(x,t) \Big] dx \, dt
  &=& 0 \nonumber
\end{eqnarray}
is satisfied for all $(v_H,q_h) \in X_H \times Y_h$. The related error
estimate now follows from the abstract estimate
\eqref{eq:error-estimate-XH}, i.e., with $c_2^B=\sqrt{2}$
and $\widetilde{c}_S=1$ this gives
\[
  \| u - \widetilde{u}_H \|_{X,h} \leq
  \big( 1 + 2\sqrt{2} \big)
  \inf\limits_{v_H \in X_H} \| u - v_H \|_X \, .
\]
Hence, assuming $u \in H^s(Q)$ for some $s \in [1,2]$ we conclude
the error estimate, e.g., \cite{Steinbach:2015},
\begin{equation*}
  \| \nabla_x (u - \widetilde{u}_H) \|_{L^2(Q)} \leq
  \| u - \widetilde{u}_H \|_{X,h} \leq c H^{s-1}|u|_{H^s(Q)}.
\end{equation*}
We also define the ansatz space $X_{\overline{H}} := Y_h \cap X$ and compute 
$(\widetilde{u}_{\overline{H}},\overline{p}_h) \in X_{\overline{H}} \times Y_h$
such that
\begin{eqnarray}
  \int_Q \nabla_x \overline{p}_h(x,t) \cdot \nabla_x q_h(x,t) \, dx \, dt
  \hspace*{6cm}
  && \nonumber \\
  +  \int_Q \Big [\partial_t \widetilde{u}_{\overline{H}}(x,t) \, q_h(x,t)  +
  \nabla_x \widetilde{u}_{\overline{H}}(x,t)
  \cdot \nabla_x q_h(x,t) \Big] dx \, dt
  &=& \langle f , q_h \rangle_Q, \label{VF heat Xbar} \\
  \int_Q \Big[ \partial_t v_{\overline{H}}(x,t) \, \overline{p}_h(x,t) +
  \nabla_x v_{\overline{H}}(x,t) \cdot \nabla_x
  \overline{p}_h(x,t) \Big] dx \, dt 
  &=& 0 \nonumber
\end{eqnarray}
is satisfied for all $(v_{\overline{H}},q_h) \in
X_{\overline{H}} \times Y_h$. As before we can prove the discrete
inf-sup stability condition
\[
  \| u_{\overline{H}} \|_{X,h} \leq \sup\limits_{0 \neq q_h \in Y_h}
  \frac{\spf{Bu_{\overline{H}}}{q_h}_Q}{\|q_h\|_Y} \quad \mbox{for all} \;
  u_{\overline{H}} \in X_{\overline{H}} ,
\]
which ensures unique solvability of \eqref{VF heat Xbar}, and we write
the abstract saturation condition \eqref{saturation abstract} as
\[
  \| u - \widetilde{u}_{\overline{H}} \|_{X,h} \leq \eta \,
  \| u - \widetilde{u}_H \|_{X,h} \quad \mbox{for some} \; \eta \in (0,1).
\]
As in the case of the Poisson equation we now use global error estimator
\[
  \eta_h^2 = \| \nabla_x p_h \|^2_{L^2(Q)} = \int_Q |\nabla_x p_h(x,t)|^2
  \, dx \, dt = \sum\limits_{\ell=1}^{N_h} \int_{q_\ell^h}
  |\nabla_x p_h(x,t)|^2 \, dx \, dt =
  \sum\limits_{\ell=1}^{N_h} \eta_{h,\ell}^2
\]
with the local error indicators
\[
  \eta_{h,\ell}^2 = \int_{q_\ell^h} |\nabla_x p_h(x,t)|^2 \, dx \, dt, \quad
  \eta_{H,j}^2 = \sum\limits_{\tau_\ell^h \subset \tau_j^H} \eta_{h,\ell}^2
\]
to drive an adaptive refinement algorithm, using D\"orfler marking with
$\theta = 0.5$ and newest vertex bisection,
\cite{Doerfler:1996}, and we compare these results
with those obtained by a uniform refinement strategy.
We consider the Galerkin variational formulation
\eqref{eq:DBVP Parabolic FEM} for the ansatz space
$X_H = S_H^1(Q) \cap X$ of piecewise linear and continuous
basis functions, and the test space $Y_H = S_H^2(Q) \cap Y$ of piecewise
second order and continuous basis functions. The latter corresponds to
the test space $Y_h = S_h^1(Q) \cap Y$ of piecewise linear continuous
basis functions, which are defined with respect to a refined mesh with
local mesh size $h=\frac{1}{2}H$.
As before, all computations were done in Netgen/NGSolve 
\cite{Schoeberl:2014} where
we used the sparse direct solver package Umfpack to solve the
resulting linear system.

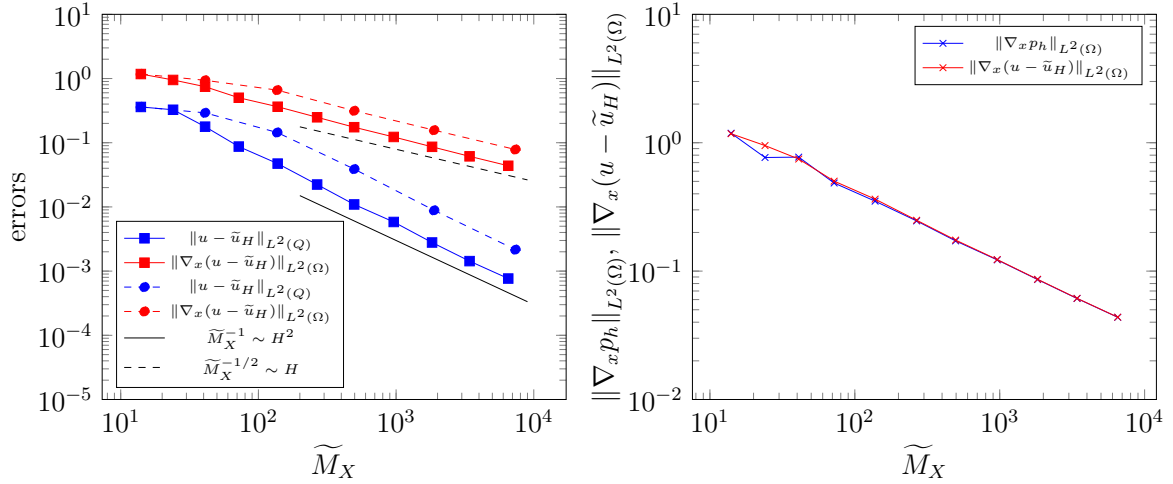
\begin{figure}[!htbp]
	\centering
	\begin{subfigure}[b]{0.45\textwidth}
		\begin{tikzpicture}[scale = 0.9]
			\begin{axis}[
				xmode = log,
				ymode = log,
				xlabel= $\widetilde{M}_X$,
				ylabel= errors,
				ymin = 0.00001,
				ymax = 10,
				legend style={font=\tiny}, legend pos = south west]
				\addplot[mark = square*,blue] table [col sep=&, y=errL2_u, x=nv]{tables/SimData_LSHeat_ex1_adaptive.dat};
				\addlegendentry{$\norm{u - \widetilde{u}_H}_{L^2(Q)}$}
			
				\addplot[mark = square*,red] table [col sep=&, y=errH1_u, x=nv]{tables/SimData_LSHeat_ex1_adaptive.dat};
				\addlegendentry{$\norm{\nabla_x (u-\widetilde{u}_H)}_{L^2(\Omega)}$}
				
				\addplot[mark = *,blue, dashed] table [col sep=&, y=errL2_u, x=nv]{tables/SimData_LSHeat_ex1_uniform.dat};
				\addlegendentry{$\norm{u - \widetilde{u}_H}_{L^2(Q)}$}
				
				\addplot[mark = *,red, dashed] table [col sep=&, y=errH1_u, x=nv]{tables/SimData_LSHeat_ex1_uniform.dat};
				\addlegendentry{$\norm{\nabla_x (u-\widetilde{u}_H)}_{L^2(\Omega)}$}
				
				\addplot[
				domain = 200:9000,
				samples = 10,
				thin,
				black,
				] {3*x^(-1)};
				\addlegendentry{$\widetilde{M}_X^{-1}\sim H^{2}$}
				
				\addplot[
				domain = 200:9000,
				samples = 10,
				thin,
				black,
				dashed,
			] {2.5*x^(-1/2)};
			\addlegendentry{$\widetilde{M}_X^{-1/2}\sim H$}
			\end{axis}
		\end{tikzpicture}
		\caption{Errors $\| \nabla_x (u - \widetilde{u}_H) \|_{L^2(Q)}$ and
          $\| u - \widetilde{u}_H \|_{L^2(Q)}$ for uniform and adaptive
          refinement strategies.}
		\label{abb:ex_1_error_plot}
	\end{subfigure}
	~
	\begin{subfigure}[b]{0.45\textwidth}
		\begin{tikzpicture}[scale = 0.9]
			\begin{axis}[
				xmode = log,
				ymode = log,
				xlabel= $\widetilde{M}_X$,
				ylabel= {$\norm{\nabla_x p_h}_{L^2(\Omega)}$, $\norm{\nabla_x (u-\widetilde{u}_H)}_{L^2(\Omega)}$},
				ymin = 0.01,
				ymax = 10,
				legend style={font=\tiny}, legend pos = north east]
				\addplot[mark = x,blue] table [col sep=&, y=L2_dxph, x=nv]{tables/SimData_LSHeat_ex1_adaptive.dat};
				\addlegendentry{$\norm{\nabla_x p_h}_{L^2(\Omega)}$}
			
				\addplot[mark = x,red] table [col sep=&, y=errH1_u, x=nv]{tables/SimData_LSHeat_ex1_adaptive.dat};
				\addlegendentry{$\norm{\nabla_x (u-\widetilde{u}_H)}_{L^2(\Omega)}$}
			\end{axis}
		\end{tikzpicture}
		\caption{History of the norms $\norm{\nabla_x p_h}_{L^2(Q)}$, $\norm{\nabla_x (u-
			\widetilde{u}_H)}_{L^2(Q)}$ for an adaptive refinement strategy.}
		\label{abb:ex_1_plot_ph_u-uh}
	\end{subfigure}
	
	\caption{Convergence results in the case of a smooth solution for the
        heat equation.}%
\end{figure}

\begin{figure}[!htbp]
  \centering
  \begin{subfigure}[b]{0.31\textwidth}
    \centering
    \includegraphics[width = 0.8\textwidth]{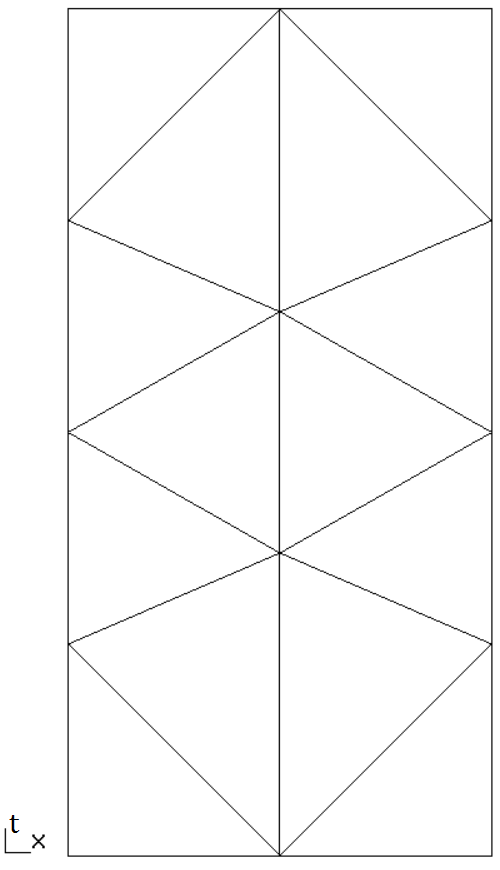}
    \caption{\footnotesize Initial mesh}\label{abb:ex1_start_mesh}
   \end{subfigure}
   ~
   \begin{subfigure}[b]{0.31\textwidth}
     \centering
     \includegraphics[width = 0.8\textwidth]{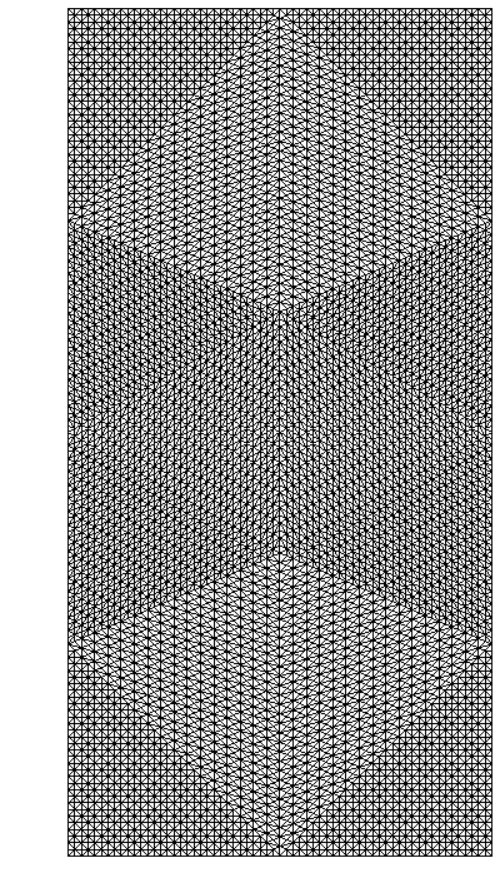}
     \caption{\footnotesize
       Uniform refinement~$L=5$}\label{abb:ex1_uniform_mesh_L5}
   \end{subfigure}
   ~
   \begin{subfigure}[b]{0.31\textwidth}
     \centering
     \includegraphics[width = 0.8\textwidth]{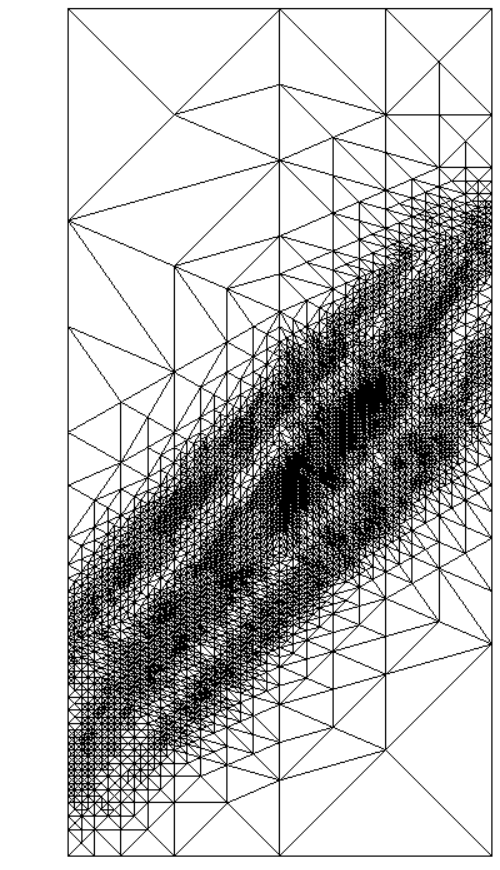}
     \caption{\footnotesize Adaptive refinement $L=10$}
   \end{subfigure}
   \caption{Space-time finite element meshes.}
	\label{abb:heat_ex1_meshes}
\end{figure}

In the first example we use the one-dimensional spatial domain
$\Omega=(0,3)$ and the time horizon $T=6$, i.e., we have the
space-time domain $ Q:= (0,3) \times (0,6) \subset {\mathbb{R}}^2$.
As exact solution we consider the smooth function
\begin{equation}\label{eq:example-smooth}
  u(x,t) :=
  \begin{cases}
    \frac{1}{2} (t-x-2)^3 (x-t)^3 \sin \frac{\pi}{3} x
    & \mbox{for} \; x \leq t~\text{and }t-x\leq 2,\\[1mm]
    0 & \text{else},
  \end{cases}
\end{equation}
and we compute $f=\partial_tu-\Delta_xu$ accordingly.
Since $u$ is smooth we expect optimal orders of convergence, i.e.,
$\mathcal{O}(H^2)$ when measuring the error in $L^2(Q)$, and
$\mathcal{O}(H)$ in the energy norm $\| \nabla_x \cdot \|_{L^2(Q)}$.
These rates are confirmed by the numerical results for both a uniform
and an adaptive refinement strategy, as shown in
Fig.~\ref{abb:ex_1_error_plot}.
In Fig.~\ref{abb:ex_1_plot_ph_u-uh} we present a comparison of the
error $\| \nabla_x(u - \widetilde{u}_H) \|_{L^2(Q)}$ with the error
estimator $\| \nabla_x p_h \|_{L^2(Q)}$, where the curves are almost
parallel. Moreover, in Table \ref{tab:ex1_comparison_of_errors} we
provide a comparison of the errors $\| \nabla_x (u - \widetilde{u}_H) \|_{L^2(Q)}$
for both the uniform and adaptive refinement strategy. We observe that
in the adaptive case less degrees of freedom are needed to
reach the same level of accuracy than in the case of a uniform refinement.
Finally, in Fig.~\ref{abb:heat_ex1_meshes} we present the related finite
element meshes. Note that the solution is smooth but behaves similar
as a wave. This motivates to use a mesh which is adaptive in space and time.

%\begin{figure}[!http]
  %\centering
  %\begin{subfigure}[b]{0.45\textwidth}
    %\includegraphics[scale = 0.5]{figures/heat_ex1_errors.pdf}%
%
    %\caption{Errors $\| u - \widetilde{u}_H \|_Y$ and
          %$\| u - \widetilde{u}_H \|_{L^2(Q)}$ for uniform and adaptive
          %refinement strategies}
	%\label{abb:ex_1_error_plot}
	%\end{subfigure}
	%~
	%\begin{subfigure}[b]{0.45\textwidth}
          %\includegraphics[scale = 0.5]{figures/heat_ex1_estimatorVSerror.pdf}
          %\caption{History of the norms $\norm{p_h}_Y$, $\norm{u-
	%\widetilde{u}_H}_Y$ for an adaptive refinement strategy}
	%\label{abb:ex_1_plot_ph_u-uh}
	%\end{subfigure}
	%\caption{Convergence results in the case of a smooth solution for the
        %heat equation}
%\end{figure}

\begin{table}[!htbp]
  \begin{center}
    \begin{tabular}{rrrrrr} \toprule
\multicolumn{3}{c}{uniform} & \multicolumn{3}{c}{adaptive} \\
\cmidrule{2-3} \cmidrule{5-6} 
$L$ & $\widetilde{M}_X$ & $\norm{\nabla_x (u-\widetilde{u}_H)}_{L^2(Q)}$ & $L$ & $\widetilde{M}_X$ & $\norm{\nabla_x (u-\widetilde{u}_H)}_{L^2(Q)}$\\
\midrule 
0 & 14 & 1.175e+00 &  0 & 14 & 1.175e+00 \\ 
1 & 41 & 9.391e-01 &  2 & 41 & 7.463e-01 \\ 
2 & 137 & 6.597e-01 & 4 & 138 & 3.637e-01 \\
3 & 497 & 3.144e-01 & 6 & 496 & 1.745e-01 \\
4 & 1889 & 1.568e-01 & 8 & 1825 & 8.636e-02 \\
5 & 7361 & 7.840e-02 & 10 & 6524 & 4.377e-02 \\ 
\bottomrule
\end{tabular}

  \end{center}
  \caption{Comparison of the errors between uniform and adaptive refinement.}
  \label{tab:ex1_comparison_of_errors}
\end{table}

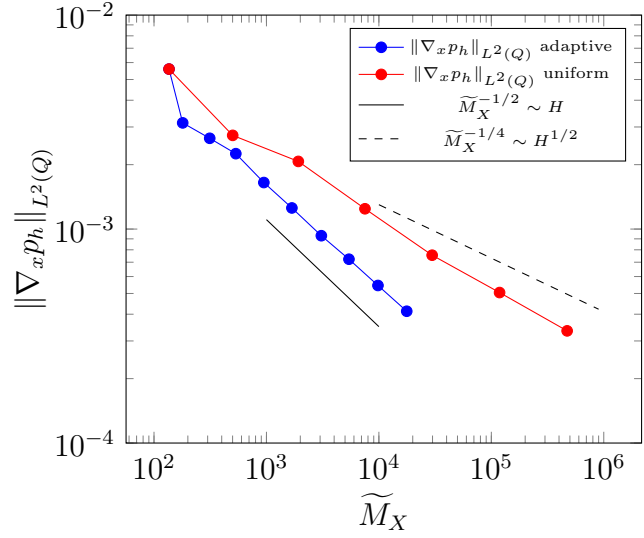
\begin{figure}[!htbp]
	\centering
	\begin{tikzpicture}
		\begin{axis}[
			xmode = log,
			ymode = log,
			xlabel= $\widetilde{M}_X$,
			ylabel= {$\norm{\nabla_x p_h}_{L^2(Q)}$},
			ymin = 0.0001,
			ymax = 0.01,
			legend style={font=\tiny}, legend pos = north east]
			\addplot[mark = *,blue] table [col sep=&, y=L2_dxph, x=nv]{tables/SimData_LSHeat_ex2_adaptive.dat};
			\addlegendentry{$\| \nabla_x p_h\|_{L^2(Q)}$ adaptive}
						
			\addplot[mark = *,red] table [col sep=&, y=L2_dxph, x=nv]{tables/SimData_LSHeat_ex2_uniform.dat};
			\addlegendentry{$\| \nabla_x p_h\|_{L^2(Q)}$ uniform}
						
			\addplot[
				domain = 1000:10000,
				samples = 10,
				thin,
				black,
			] {0.035*x^(-1/2)};
			\addlegendentry{$\widetilde{M}_X^{-1/2}\sim H$}
						
			\addplot[
				domain = 10000:900000,
				samples = 10,
				thin,
				black,
				dashed,
			] {0.013*x^(-1/4)};
			\addlegendentry{$\widetilde{M}_X^{-1/4}\sim H^{1/2}$}
		\end{axis}
	\end{tikzpicture}
	\caption{Error estimator $\| \nabla_x p_h \|_{L^2(Q)}$ in the case of a discontinuous
					 right hand side.}%
	\label{abb:heat_ex2_estimators}%
\end{figure}

\begin{figure}[!htbp]
  \centering
  \begin{subfigure}[b]{0.45\textwidth}
    \includegraphics[width = 0.8\textwidth]{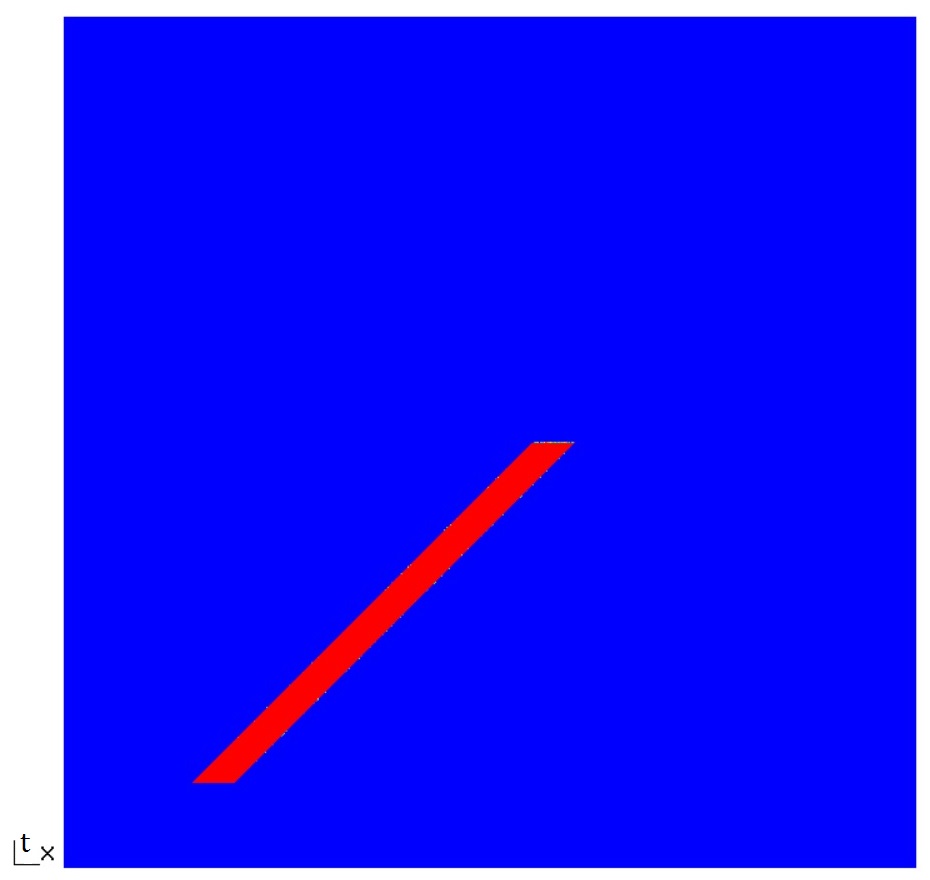}%
    \caption{Discontinuous right hand side \eqref{heat discontinuous f}}
    \label{abb:heat_ex2_f}
  \end{subfigure}
  ~
  \begin{subfigure}[b]{0.45\textwidth}
    \includegraphics[width = 0.8\textwidth]{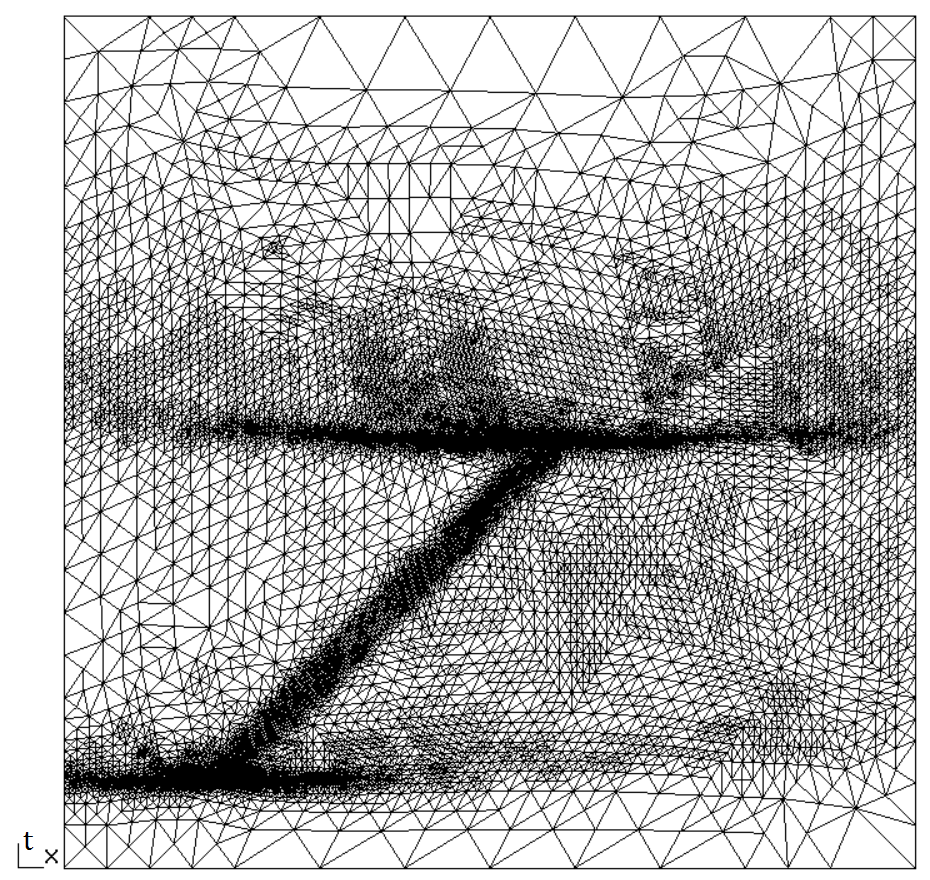}
    \caption{Adaptive mesh with $\widetilde{M}_X=17689$}
    \label{abb:heat_ex2_adaptivemesh}
  \end{subfigure}
  \caption{Singular solution of the heat equation for discontinuous right hand
    side.}
\end{figure}

As a second example we consider the heat equation \eqref{eq:Parabolic DBVP}
in the space-time domain
$Q=(0,1)^2 \subset {\mathbb{R}}^2$, i.e., $\Omega=(0,1)$ and $T=1$ with
homogeneous Dirichlet and initial conditions, but a discontinuous right
hand side
\begin{equation}\label{heat discontinuous f}
  f(x,t) =
  \begin{cases}
    1 \quad (x,t) \in \left\{ (x,t) \in (0,1)\times
      \left(\frac{1}{10},\frac{1}{2}\right)
      :~x-\frac{1}{10} \leq t \leq x-\frac{1}{20}\right\},\\
    0 \quad \text{else},
  \end{cases}
\end{equation}
as depicted in Fig.~\ref{abb:heat_ex2_f}.
Note that a similar example was also considered in
\cite{FuehrerKarkulik:2021}. Since the exact solution is unknown,
in Fig.~\ref{abb:heat_ex2_estimators} we provide the results for the
error estimator $\| \nabla_x p_h \|_{L^2(Q)}$ which is equivalent to the
error $\| \nabla_x (u - \widetilde{u}_H)\|_{L^2(Q)}$. As already observed in
\cite[Section 5.2.3]{FuehrerKarkulik:2021}, in the case of uniform
refinement we obtain a reduced rate of $\widetilde {M}_X^{-\frac{1}{4}}$
which corresponds to $\mathcal{O}(H^\frac{1}{2})$, while in the case
of an adaptive refinement we have the optimal rate of $\widetilde {M}_X^{-\frac{1}{2}}$,
i.e., $\mathcal{O}(H)$. 
In Fig.~\ref{abb:heat_ex2_adaptivemesh} we provide the adaptive mesh at
level $L=9$.

%\begin{figure}[!htbp]
	%\centering
	%\begin{subfigure}[b]{0.31\textwidth}
		%\centering
		%\begin{tikzpicture}[scale = 0.5]
		%\begin{axis}[
			%xmode = log,
			%ymode = log,
			%xlabel=dofs (coarse mesh),
			%ylabel= {$\norm{\nabla_x p_h}_{L^2(Q)}$},
			%ymin = 0.0001,
			%ymax = 0.01,
			%legend style={font=\tiny}, legend pos = north east]
			%\addplot[mark = *,blue] table [col sep=&, y=L2_dxph, x=nv]{tables/SimData_LSHeat_ex2_adaptive.dat};
			%\addlegendentry{$\| \nabla_x p_h\|_{L^2(Q)}$ adaptive}
						%
			%\addplot[mark = *,red] table [col sep=&, y=L2_dxph, x=nv]{tables/SimData_LSHeat_ex2_uniform.dat};
			%\addlegendentry{$\| \nabla_x p_h\|_{L^2(Q)}$ uniform}
						%
			%\addplot[
				%domain = 1000:10000,
				%samples = 10,
				%thin,
				%black,
			%] {0.035*x^(-1/2)};
			%\addlegendentry{$N^{-1/2}$}
						%
			%\addplot[
				%domain = 10000:900000,
				%samples = 10,
				%thin,
				%black,
				%dashed,
			%] {0.013*x^(-1/4)};
			%\addlegendentry{$N^{-1/4}$}
		%\end{axis}
		%\end{tikzpicture}
		%\caption{Error estimator $\| \nabla_x p_h \|_{L^2(Q)}$ in the case of a discontinuous
					 %right hand side}%
		%\label{abb:heat_ex2_estimators}
	%\end{subfigure}
	%~
	%\begin{subfigure}[b]{0.31\textwidth}
		%\includegraphics[width = 0.6\textwidth]{figures/LS_heat_ex2_f.jpg}%
    %\caption{Discontinuous right hand side \eqref{heat discontinuous f}}
    %\label{abb:heat_ex2_f}
	%\end{subfigure}
	%~
	%\begin{subfigure}[b]{0.31\textwidth}
		%\includegraphics[width = 0.6\textwidth]{figures/LS_heat_ex2_adaptivemesh_L9.png}
    %\caption{Adaptive mesh with $17689$ dofs}
    %\label{abb:heat_ex2_adaptivemesh}
	%\end{subfigure}
%
%\end{figure}

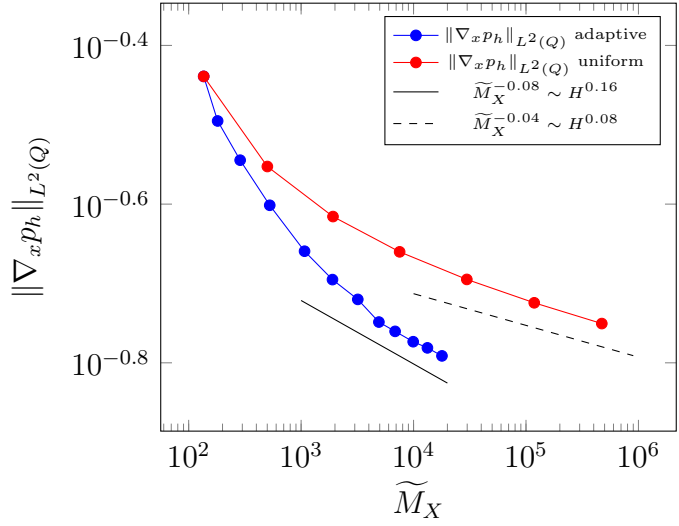
\begin{figure}
	\centering
	\begin{tikzpicture}
		\begin{axis}[
			xmode = log,
			ymode = log,
			xlabel= $\widetilde{M}_X$,
			ylabel = {$\norm{\nabla_x p_h}_{L^2(Q)}$},
			ylabel near ticks,
			ymin = 0.13,
			ymax = 0.45,
			legend style={font=\tiny}, legend pos = north east]
			\addplot[mark = *,blue] table [col sep=&, y=L2_dxph, x=nv]{tables/SimData_LSHeat_ex3_adaptive.dat};
			\addlegendentry{$\| \nabla_x p_h\|_{L^2(Q)}$ adaptive}
						
			\addplot[mark = *,red] table [col sep=&, y=L2_dxph, x=nv]{tables/SimData_LSHeat_ex3_uniform.dat};
			\addlegendentry{$\| \nabla_x p_h\|_{L^2(Q)}$ uniform}
			
			\addplot[
				domain = 1000:20000,
				samples = 10,
				thin,
				black,
			] {0.33*x^(-0.08)};
			\addlegendentry{$\widetilde{M}_X^{-0.08}\sim H^{0.16}$}
			
			\addplot[
				domain = 10000:900000,
				samples = 10,
				thin,
				black,
				dashed,
			] {0.28*x^(-0.04)};
			\addlegendentry{$\widetilde{M}_X^{-0.04}\sim H^{0.08}$}
		\end{axis}
	\end{tikzpicture}
	\caption{Error estimator $\| \nabla_x p_h \|_{L^2(Q)}$ in the case of non compatible
    initial conditions $u_0 \not\in H^1_0(\Omega)$.}%
  \label{abb:heat_ex3_errors}%
\end{figure}

As a last example for the heat equation we consider again the
space-time domain $Q = (0,1)^2$, with $f(x,t)=2$ for $(x,t)\in Q$,
and $u_0(x)=1$ for $x \in (0,1)$. Obviously, we have
$u_0 \in L^2(0,1)$, but $u_0 \not\in H^1_0(\Omega)$, i.e., there is
no compatibility with the homogeneous Dirichlet boundary conditions at
$t = 0$. This results in a reduced order of convergence. In this example,
we choose the parameter $\theta = 0.9$ in the D\"orfler marking strategy.
Since the exact solution is not known, in Fig.~\ref{abb:heat_ex3_errors}
we provide the results for the error indicator
$\| \nabla_x p_h \|_{L^2(Q)}$ which is equivalent to the error
$\| \nabla_x (u - \widetilde{u}_H) \|_{L^2(Q)}$. We observe a rate of $\mathcal{O}(H^{0.08})$
in the case of a uniform refinement, and of $\mathcal{O}(H^{0.16})$ for the
adaptive refinement strategy. These rates coincide with those observed in 
\cite[Section 5.2.4]{FuehrerKarkulik:2021}.

\begin{remark}
  Similar as in the case for the Poisson equation we may also define
  \[
    X := \Big \{ v \in L^2(0,T;L^2(\Omega)) \cap H^1_{0,}(0,T;H^{-1}(\Omega)) :
    \partial_t v - \Delta_x v \in L^2(Q) \Big \}, \quad Y = L^2(Q).
  \]
  In this case we have $B = \partial_t - \Delta_x : X \to L^2(Q)$, and
  $A := I : L^2(Q) \to L^2(Q)$. For a conforming space-time finite element
  discretization one can use tensor-product meshes with quadratic B splines
  in the spatial directions, and piecewise linear continuous basis
  functions in the temporal direction.
\end{remark}

\section{Hyperbolic Dirichlet boundary value problem}
In this section we apply the abstract theory to the 
Dirichlet boundary value problem for the wave equation,
\begin{subequations}\label{eq:Hyperbolic DBVP}
  \begin{alignat}{2}
    \Box u(x,t) \, := \, \partial_{tt} u(x,t) - \Delta_x u(x,t)
    & = f(x,t) & \quad & \text{for } (x,t) \in Q := \Omega \times (0,T), \\
    u(x,t)
    & = 0 & \quad & \text{for } (x,t) \in \Sigma := \Gamma \times (0,T),\\
    u(x,0) = \partial_t u(x,t)_{|t=0} & = 0 & \quad & \text{for } x \in
    \Omega,
  \end{alignat}
\end{subequations}
where, as in the parabolic case, $\Omega \subset \R^n$, $n = 1,2,3$ is a
bounded Lipschitz domain with boundary $\Gamma:= \partial \Omega$,
and $T>0$ is a given time horizon.
Following \cite{SteinbachZank:2020}, the space-time variational
formulation of the wave equation \eqref{eq:Hyperbolic DBVP} is to find
$u \in H^{1,1}_{0;0,}(Q)$ such that
\begin{equation}\label{VF wave}
  - \langle \partial_t u , \partial_t q \rangle_{L^2(Q)} +
  \langle \nabla_x u , \nabla_x q \rangle_{L^2(Q)} =
  \langle f , q \rangle_Q
\end{equation}
is satisfied for all $q \in H^{1,1}_{0;,0}(Q)$. Here we use the
anisotropic Sobolev space
\[
H^{1,1}_{0;0,}(Q) := L^2(0,T;H^1_0(\Omega)) \cap H^1_{0,}(0,T;L^2(\Omega)),
\]
where $H^1_{0,}(0,T;L^2(\Omega))$ covers the zero initial condition
$u(x,0)=0$ for $x \in \Omega$, while $L^2(0,T;H^1_0(\Omega))$ includes
the homogeneous Dirichlet boundary condition on $\Sigma$. Note that the
second initial condition $\partial_t u(x,t)_{|t=0}=0$ for $x \in \Omega$
enters the variational formulation \eqref{VF wave} in
a natural way. A norm in $H^{1,1}_{0;0,}(Q)$ is given by the graph norm
\[
  \| u \|_{H^{1,1}_{0;0,}(Q)} := \sqrt{ \| \partial_t u \|_{L^2(Q)}^2 +
  \| \nabla_x u \|^2_{L^2(Q)}} \, = \, |u|_{H^1(Q)} .
\]
Note that $H^{1,1}_{0;,0}(Q)$ is defined accordingly, but with a zero
terminal condition $q(x,T)=0$ for $x \in \Omega$. Although the right
hand side of the variational formulation is well defined for
$ f \in [H^{1,1}_{0;,0}(Q)]^*$, in order to ensure a unique solution
$ u \in H^{1,1}_{0;0,}(Q)$ we have to assume $f \in L^2(Q)$, e.g.,
\cite{Ladyzhenskaya:1985,SteinbachZank:2020}. In fact, the
solution operator mapping $f \in L^2(Q)$ to the solution
$u \in H^{1,1}_{0;0,}(Q)$ of the variational formulation
\eqref{VF wave} does not define an isomorphism. Instead we have to
enlarge the ansatz space in order to
incorporate the second initial condition
$\partial_t u(x,t)_{|t=0}=0$ in an appropriate way.
In what follows we will consider a generalized variational formulation
of the wave equation, see \cite{SteinbachZank:2022}: For the
enlarged space-time domain $Q_-:=\Omega \times (-T,T)$, and for
$u\in L^2(Q)$ we define the zero extension
\[
  \widetilde{u}(x,t) :=
  \begin{cases}
	u(x,t) & \mbox{for} \; (x,t)\in Q, \\
	0, & \mbox{else.}
      \end{cases}
\]
The application of the wave operator $\Box \widetilde{u}$ on $Q_-$ will
be formulated as a distribution, i.e., for
$\varphi \in C_0^\infty(Q_-)$ we define
\[
  \langle \Box \widetilde{u}, \varphi \rangle_{Q_-}
  :=
  \int_{Q_-} \widetilde{u}(x,t) \, \Box \varphi(x,t) \, dx \, dt =
  \int_Q u(x,t) \, \Box \varphi(x,t) \, dx \, dt.  
\]
Now we are in the position to introduce the space 
\[
  \mathcal{H}(Q) := \Big \{
  u=\widetilde{u}_{|_{Q}} : \widetilde{u} \in L^2(Q_-) , \,
  \widetilde{u}_{|\Omega \times (-T,0)}=0, \,
  \Box \widetilde{u} \in [H_0^1(Q_-)]^* \Big \},
\]
with the graph norm
\[
  \|u\|_{\mathcal{H}(Q)} :=
  \sqrt{\|u\|_{L^2(Q)}^2 + \|\Box \widetilde{u} \|_{[H_0^1(Q_-)]^*}^2} \, .
\]
The normed vector space $(\mathcal{H}(Q),\|\cdot\|_{\mathcal{H}(Q)})$ is a
Banach space, and it holds true that, see \cite[Lemma 3.5]{SteinbachZank:2022},
$H_{0;0,}^{1,1}(Q)\subset \mathcal{H}(Q)$ i.e.,
\begin{equation}\label{norm X H1}
  \| \Box \widetilde{u} \|_{[H^1_0(Q_-)]^*} \leq \| u \|_{H^{1,1}_{0;0,}(Q)}
  \quad \mbox{for all} \; u \in H^{1,1}_{0;0,}(Q).
\end{equation}
Therefore, we can consider the space
\[
  \mathcal{H}_{0;0,}(Q) :=
  \overline{H_{0;0,}^{1,1}(Q)}^{\| \cdot \|_{\mathcal{H}(Q)}} \subset
  \mathcal{H}(Q)
\]
which will serve as ansatz space. For $u \in {\mathcal{H}}_{0;0,}(Q)$,
an equivalent norm is given as, see \cite[Lemma 3.6]{SteinbachZank:2022},
\[
  \| u \|_{{\mathcal{H}}_{0;0,}(Q)} = \| \Box \widetilde{u}
  \|_{[H^1_0(Q_-)]^*} .
\]
For given $f \in [H_{0;,0}^{1,1}(Q)]^*$ we now
consider the variational formulation
to find $u \in \mathcal{H}_{0;0,}(Q)$ such that 
\begin{equation}\label{eq:generalized-VF-wave}
  b(u,q) := \langle B u , q \rangle_Q =
  \langle \Box \widetilde{u} , {\mathcal{E}} q \rangle_{Q_-}
  = \langle f , q \rangle_{Q} \quad
  \mbox{for all}  \; q\in H_{0;,0}^{1,1}(Q),
\end{equation}
where $\mathcal{E} : H_{0;,0}^{1,1}(Q) \to H_0^1(Q_-)$ is a suitable
extension operator, e.g., reflection in time. Note that we have, see
  {\rm \cite[Lemma 3.5]{SteinbachZank:2022}},
  \[
    \langle \Box \widetilde{u} , {\mathcal{E}} q \rangle_{Q_-}
    =
    - \langle \partial_t u , \partial_t q \rangle_{L^2(Q)} +
    \langle \nabla_x u , \nabla_x q \rangle_{L^2(Q)}, \quad
    u \in H^{1,1}_{0;0,}(Q) \subset {\mathcal{H}}_{0;0,}(Q),
    q \in H^{1,1}_{0;,0}(Q) .
  \]
We conclude that the bilinear form within the variational
formulation \eqref{eq:generalized-VF-wave} is bounded for
all $u \in {\mathcal{H}}_{0;0,}(Q)$ and $q \in H^{1,1}_{0;,0}(Q)$,
and satisfies the inf-sup stability condition
\begin{equation}\label{inf-sup generalized}
  \| u \|_{\mathcal{H}_{0;0,}(Q)} \, = \, 
  \sup_{0\neq q\in H_{0;,0}^{1,1}(Q)}
  \frac{\langle \Box \widetilde{u} , {\mathcal{E}} q \rangle_{Q_-}}
  {\| q \|_{H_{0;,0}^{1,1}(Q)}} \quad  \mbox{for all} \;
  u \in \mathcal{H}_{0;0,}(Q) .
\end{equation}
In view of the abstract setting we therefore have
\[
  X := \mathcal{H}_{0;0,}(Q), \quad
  Y := H_{0;,0}^{1,1}(Q), \quad H:=L^2(Q), \quad B : X \to Y^* .
\]
Finally, we define
\[
  \langle A p , q \rangle_Q :=
  \langle \partial_t p , \partial_t q \rangle_{L^2(Q)} +
  \langle \nabla_x p , \nabla_x q \rangle_{L^2(Q)} \quad
  \mbox{for all} \; p,q \in Y = H^{1,1}_{0;,0}(Q),
\]
which corresponds to the space-time Laplacian. Thus we can apply the
abstract theory as given in Section~\ref{sec:abstract setting}.
The variational formulation \eqref{VF Coupled system Su=g} now reads
to find $(u,p) \in X \times Y$ such that
\begin{equation}\label{system wave}
  \langle \partial_t p , \partial_t q \rangle_{L^2(Q)} +
  \langle \nabla_x p , \nabla_x q \rangle_{L^2(Q)} +
  \langle \Box \widetilde{u} , {\mathcal{E}} q \rangle_{Q_-} =
  \langle f , q \rangle_Q , \quad
  \langle \Box \widetilde{v} , {\mathcal{E}} p \rangle_{Q_-} = 0
\end{equation}
is satisfied for all $(v,q) \in X \times Y$.

Let $X_H \subset H^{1,1}_{0;0,}(Q) \subset {\mathcal{H}}_{0;0,}(Q) = X$ be
the conforming finite element space of piecewise linear and continuous 
basis functions which are defined with respect to an admissible locally
quasi-uniform decomposition of the space-time domain $Q$ into
shape-regular simplicial finite elements of mesh size $H$. Moreover,
let $Y_h \subset H^{1,1}_{0;,0}(Q)$ be a second finite element space of
piecewise linear and continuous basis functions which are defined
with respect to a refined decomposition of the space-time domain
into finite elements of mesh size $h$. The Galerkin discretization of
\eqref{system wave} then reads to find $(u_H,p_h) \in X_H \times Y_h$
such that
\begin{equation}\label{wave FEM 1}
  \langle \partial_t p_h , \partial_t q_h \rangle_{L^2(Q)} +
  \langle \nabla_x p_h , \nabla_x q_h \rangle_{L^2(Q)} -
  \langle \partial_t u_H , \partial_t q_h \rangle_{L^2(Q)} +
  \langle \nabla_x u_H , \nabla_x q_h \rangle_{L^2(Q)} =
  \langle f , q_h \rangle_Q
\end{equation}
and
\begin{equation}\label{wave FEM 2}
- \langle \partial_t v_H , \partial_t p_h \rangle_{L^2(Q)} +
  \langle \nabla_x v_H , \nabla_x p_h \rangle_{L^2(Q)} = 0
\end{equation}
is satisfied for all $(v_H,q_h) \in X_H \times Y_h$.
Using the relation $\| u \|_{\mathcal{H}_{0;0,}(Q)} \leq
\| u \|_{H^{1,1}_{0;0,}(Q)} = |u|_{H^1(Q)}$ for $u\in H^{1,1}_{0;0,}(Q)$ we
conclude a best approximation result. 

\begin{theorem}
  Assume \eqref{Assumption Yh} such that the discrete inf-sup stability
  condition \eqref{inf-sup Su=g} is satisfied, and assume
  $u\in H^{1,1}_{0;0,}(Q)\cap H^s(Q)$ for some $s \in [1,2]$.
  For the unique solution
  $(u_H,p_h) \in X_H \times Y_h$ of \eqref{wave FEM 1} and \eqref{wave FEM 2} 
  there holds the error estimate
  \begin{equation*}
    \| u-u_H \|_{\mathcal{H}_{0;0,}(Q)} \leq
    \inf_{v_H\in X_H} |u-v_H|_{H^1(Q)} \leq
    c \, H^{s-1} |u|_{H^s(Q)}. 
  \end{equation*}
\end{theorem}

\noindent
If the finite element space $Y_h$ is chosen appropriately, i.e., such that
the discrete inf-sup stability condition \eqref{inf-sup Xbar} and
the saturation condition \eqref{saturation abstract} are satisfied,
$\eta = |p_h|_{H^1(Q)}$ serves as an error indicator for
$\| u-u_H \|_{\mathcal{H}_{0;0,}(Q)}$. We define the finite element
spaces $X_H = S_H^1(Q) \cap X$ and $Y_h = S_h^1(Q) \cap Y$
of piecewise linear continuous basis functions with either
$h=H/2$ or $h=H/4$, and we consider the
D\"orfler criterion \cite{Doerfler:1996} with parameter $\theta = 0.5$.
As in the previous examples we refine all marked elements
using red-green-blue refinement. 

\begin{figure}[h]
  \begin{subfigure}[b]{0.45\textwidth}
    \centering
    \begin{tikzpicture}[scale=0.6]
      \begin{axis}[
        xmode = log,
        ymode = log,
        xlabel=$M_X$ =\# DoFs (coarse mesh),
        ylabel= errors,
        ymin = 0.005,
        ymax = 5,
        legend style={font=\tiny}, legend pos = outer north east]
        \addplot[mark =*,red,dashed] table [col sep=
        &, y=errp1, x=M]{tables/tab_least-squares-uniform-h-H_2.dat};
        \addlegendentry{$|p_h|_{H^1(Q)}$ uniform}
        \addplot[mark =*,blue,dashed] table [col sep=
        &, y=erru1, x=M]{tables/tab_least-squares-uniform-h-H_2.dat};
        \addlegendentry{$|\tilde u_{H}-\overline u|_{H^1(Q)}$ uniform}
        \addplot[mark =square*,red] table [col sep=
        &, y=errp1, x=M]{tables/tab_least-squares-adaptive-h-H_2.dat};
        \addlegendentry{$|p_h|_{H^1(Q)}$ adaptive}
        \addplot[mark =square*,blue] table [col sep=
        &, y=erru1, x=M]{tables/tab_least-squares-adaptive-h-H_2.dat};

        \addlegendentry{$|\tilde u_{H}-u|_{H^1(Q)}$ adaptive}
        \addplot[domain = 300:10^5,samples = 2,dashed,thin,black]{30*x^(-1/2)};
        \addlegendentry{$M_X^{-1/2}\sim H$}
      \end{axis}
    \end{tikzpicture}
    \caption{$h=H/2$}
  \end{subfigure}
  \hfill
  \begin{subfigure}[b]{0.45\textwidth}
    \centering
    \begin{tikzpicture}[scale=0.6]
      \begin{axis}[
        xmode = log,
        ymode = log,
        xlabel=$M_X$ =\# DoFs (coarse mesh),
        ylabel= errors,
        ymin = 0.005,
        ymax = 5,
        legend style={font=\tiny}, legend pos = outer north east]
        \addplot[mark =*,red,dashed] table [col sep=
        &, y=errp1, x=M]{tables/tab_least-squares-uniform-h-H_4.dat};
        \addlegendentry{$|p_h|_{H^1(Q)}$ uniform}
        \addplot[mark =*,blue,dashed] table [col sep=
        &, y=erru1, x=M]{tables/tab_least-squares-uniform-h-H_4.dat};
        \addlegendentry{$|\tilde u_{H}-\overline u|_{H^1(Q)}$ uniform}
        \addplot[mark =square*,red] table [col sep=
        &, y=errp1, x=M]{tables/tab_least-squares-adaptive-h-H_4.dat};
        \addlegendentry{$|p_h|_{H^1(Q)}$ adaptive}
        \addplot[mark =square*,blue] table [col sep=
        &, y=erru1, x=M]{tables/tab_least-squares-adaptive-h-H_4.dat};
        \addlegendentry{$|\tilde u_{H}-u|_{H^1(Q)}$ adaptive}
        \addplot[domain = 300:10^5,samples = 2,dashed,thin,black] {30*x^(-1/2)};
        \addlegendentry{$M_X^{-1/2}\sim H$}
      \end{axis}
    \end{tikzpicture}
    \caption{$h=H/4$}
  \end{subfigure}
  \caption{Errors $| u - \widetilde{u}_H |_{H^1(Q)}$ and estimators
    $|p_h|_{H^1(Q)}$ for adaptive and uniform refinements for different
    choices of $h$ and $H$ for the smooth solution \eqref{eq:example-smooth}. }
\label{abb:wave_ex1_errors}
\end{figure}
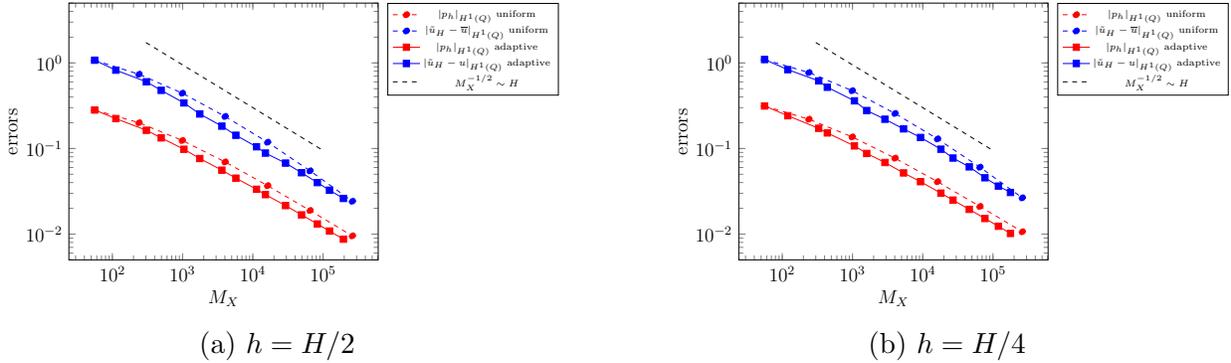

\begin{figure}[b]
  \begin{subfigure}[b]{0.29\textwidth}
    \centering 
    \includegraphics[width=\textwidth]{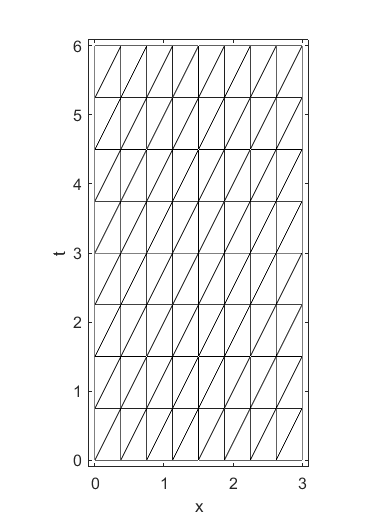}
    \caption{initial, $M_X=56$}
  \end{subfigure}
  \hfill
  \begin{subfigure}[b]{0.28\textwidth}
    \centering 
    \includegraphics[width=\textwidth]
    {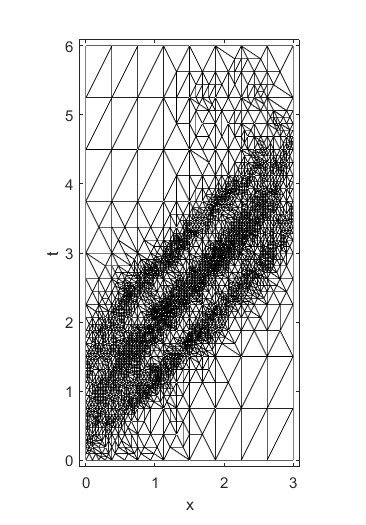}
    \caption{$h=H/2$, $M_X=3624$}
  \end{subfigure}
  \hfill
  \begin{subfigure}[b]{0.29\textwidth}
    \centering 
    \includegraphics[width=\textwidth]
    {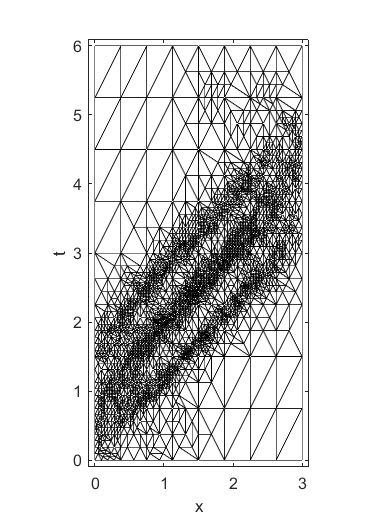}
    \caption{$h=H/4$, $M_X=2893$}
  \end{subfigure}
  \caption{Meshes from the adaptive refinement for the solution
    \eqref{eq:example-smooth}.}
  \label{abb:wave_ex1_meshes}
\end{figure}

As a first example we consider the smooth solution
\eqref{eq:example-smooth} as for the heat equation, and we compute
$f = \partial_{tt} u-\Delta_x u$. In this case we expect to see optimal
orders of convergence for the error in the energy norm.
In Fig.~\ref{abb:wave_ex1_errors} we present the numerical results
for both a uniform and an adaptive refinement, and for $h=H/2$ as
well as for $h=H/4$. In all cases we observe a linear convergence,
for both the error $|u-\widetilde{u}_H|_{H^1(Q)}$, and the
estimator $|p_h|_{H^1(Q)}$. Finally, in Fig. \ref{abb:wave_ex1_meshes}
we present the initial mesh, and two adaptively refined meshes in
the space-time domain.

\begin{figure}[!htbp]
  \centering
  \begin{tikzpicture}
    \begin{axis}[
      xmode = log,
      ymode = log,
      xlabel= $\widetilde{M}_X$,
      ylabel= {$|p_h|_{H^1(Q)}$},
      ylabel near ticks,
      ymin = 0.6,
      ymax = 0.85,
      legend style={font=\tiny}, legend pos = south west]
      \addplot[mark = *,red] table [col sep=
      &, y=L2_dxtph, x=nv]{tables/SimData_LSWave_ex2_adaptive_P1P2.dat};
      \addlegendentry{$|p_h|_{H^1(Q)}$ adaptive}
      \addplot[mark = *,blue] table [col sep=
      &, y=L2_dxtph, x=nv]{tables/SimData_LSWave_ex2_uniform_P1P2.dat};
      \addlegendentry{$|p_h|_{H^1(Q)}$ uniform}
      \addplot[
      domain = 5000:100000,
      samples = 10,
      thin,
      black,
      ] {1.13*x^(-0.05)};
      \addlegendentry{$\widetilde{M}_X^{-0.05}\sim H^{0.1}$}
      \addplot[
      domain = 10000:300000,
      samples = 10,
      dashed,
      thin,
      black,
      ] {x^(-0.03)};
      \addlegendentry{$\widetilde{M}_X^{-0.03}\sim H^{0.06}$}
    \end{axis}
  \end{tikzpicture}
  \caption{Error estimator for adaptive and uniform refinement in
    the case of a less regular solution violating the compatibility
    of initial and boundary conditions.}
  \label{abb:wave_ex2_errors}
\end{figure}
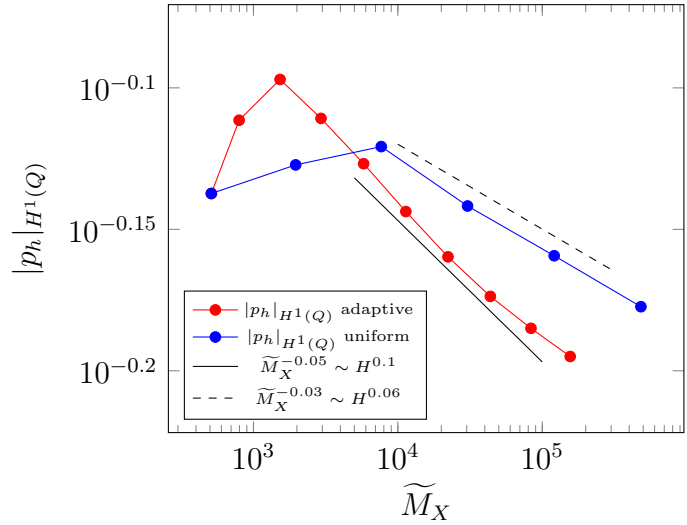

As a second example we consider the unit square $Q=(0,1)^2$ and choose
the right hand side $f$ and the inital values $u(x,0) = u_0(x)$ and
$\partial_t u(x,t)_{|t=0} = g(x)$ to be
\begin{equation*}
  f(x,t) = 2 \; \text{for } (x,t)\in Q, \quad
  u_0(x)=1 \; \text{for } x \in (0,1), \quad
  g(x) = 0 \; \text{for }x\in (0,1).
\end{equation*}
Note that the initial condition satisfies $u_0 \in L^2(0,1)$, but
$u_0 \notin H^1_0(0,1)$, i.e., there is no compatibility with the
homogeneous Dirichlet boundary conditions at $t=0$. Therefore, we
expect to see reduced orders of convergence.
This is confirmed as shown in Fig.~\ref{abb:wave_ex2_errors},
where we obeserve a rate of $\mathcal{O}(H^{0.06})$
in the uniform refinement case, using $X_H = S_H^1(Q)\cap X$ and
$Y_h = Y_H= S_H^2(Q)\cap Y$. Driving an adaptive refinement scheme with the 
D\"orfler parameter $\theta = 0.9$ the rate could be
increased to $\mathcal{O}(H^{0.1})$.

\begin{remark}
  As for the Laplace equation we can also use a boundary element
  least-squares formulation in the case of the wave equation
  \cite{Hoonhout:2023}. In particular in the one-dimensional case
  one can prove the mesh condition $H>2h$ in order to ensure
  stability, even for adaptive meshes.
\end{remark}

\section{Conclusions}
In this paper we have formulated and analyzed least-squares
methods for the numerical solution of abstract operator equations $Bu=f$.
Applications involve elliptic, parabolic, and hyperbolic problems,
with the Poisson equation, the heat equation, and the wave equation
as examples. While we assume that $B : X \to Y^*$ is an isomorphism,
and $B$ is given from the application, there is still some freedom
in the choice of the underlying spaces $X$ and $Y$. When considering
the Laplace operator $B = - \Delta$, instead of
$ B : H^1_0(\Omega) \to H^{-1}(\Omega)$ we may also consider $Y=L^2(\Omega)$,
implying $X = \{ v \in H^1_0(\Omega) : \Delta v \in L^2(\Omega)\}$, or
the other way around, i.e., an ultra weak variational formulation
using $X=L^2(\Omega)$ and
$Y = \{ v \in H^1_0(\Omega) : \Delta v \in L^2(\Omega)\}$.
It is obvious that in all of these cases we have to use appropriate
inf-sup stable finite element spaces $X_H$ and $Y_h$. But these approaches
can be applied to the heat, and the wave equation, as well, see, e.g.,
\cite{Urban:2022}.

In addition to those simple model problems can apply this methodology to
rather general problems, including flow problems, nonlinear equations, etc.
Other applications involve the least-squares formulations of boundary
integral equations \cite{Hoonhout:2023,Steinbach:LSBEM}.

In this paper we have considered the stability and error analysis of
adaptive space-time least-squares finite element methods. It is clear
that for an efficient solution of the resulting huge linear systems
of algebraic equations we need to use appropriate preconditioned
and parallel solution strategies. While the construction of these
solvers was not within the scope of this paper, this will be done in
future work. In particular, due to the structure of the linear system
to be solved, we need to have preconditioners for $A_h$, and
for the discrete Schur complement $\widetilde{S}_h = B_h^\top A_h^{-1} B_h$.
Both matrices are symmetric and positive definite, independent of the
particular choice of $B$. Possible solution strategies involve
direct methods based on factorization as used in \cite{LangerZank:2021},
or multigrid methods as considered in \cite{GanderNeumueller:2016}.
The use of efficient solution methods then also allows the numerical
solution of partial differential equations in the four dimensional
space-time domain, as we already did for optimal control problems,
e.g., \cite{2022_Langer_Steinbach_robust}.

\bigskip

\noindent
\textbf{Acknowledgement:} 
This work has been supported by the Austrian
Science Fund (FWF) under the Grant Collaborative Research Center
TRR361/F90: CREATOR Computational Electric Machine Laboratory.

\bibliographystyle{abbrv}
\bibliography{least-squares}

\end{document}